\documentstyle[12pt]{article}

\def \m{\par\medskip}
\def \b{\par\bigskip}

\newtheorem{teor}{\indent Theorem}[section]
\newtheorem{prop}{\indent Proposition}[section]
\newtheorem{lema}{\indent Lemma}[section]

\newtheorem{defin}{\indent Definition}[section]

\newtheorem{exam}{\indent Examples}[section]

\begin{document}

\title{IRREDUCIBLE REPRESENTATIONS OF THE EXCEPTIONAL CHENG-KAC SUPERALGEBRA}
\author{Consuelo Mart\'{\i}nez\thanks{Partially
supported by MTM 2010-67884-C04-01}  \\
Departamento de Matem\'aticas, Universidad de Oviedo,\\ C/ Calvo
Sotelo, s/n, 33007 Oviedo SPAIN \and Efim
Zelmanov\thanks{Partially supported the NSF}
\\Department of Mathematics,  University of California
at San Diego\\9500 Gilman Drive, La Jolla, CA 92093-0112 USA}
 \date{\it To the memory of our dear friend Hyo Chul Myung}
\maketitle


\vskip 1cm

\begin{abstract}

We classify all conformal irreducible modules of finite type over
the Cheng Kac superalgebra $CK_6$.
\end{abstract}

\section{Introduction}

The study of Lie conformal superalgebras and their representations
was initiated by V. Kac ([K2]) in view of their connections to the
free fields realizations in conformal field theory. A complete
classification of simple Lie conformal superalgebras of finite
type was achieved in [FK].  The list consists of current Lie
superalgebras, $Cur (\cal G)$, where $\cal G$ is a simple finite
dimensional Lie superalgebra; four series of Lie conformal
superalgebras of Cartan type and the exceptional Lie conformal
superalgebra $CK_6$.

\bigskip

For classification of representations of finite type of current
Lie superalgebras and Lie superalgebras of Cartan type see [BKLR],
[BKL1], [BKL2], [CK1].

In this paper we classify all conformal irreducible modules of
finite type over the superalgebra $CK_6$.  We use this classification and
the results of [MZ4] to classify conformal irreducible Jordan bimodules of finite type
over the Jordan superalgebra  $JCK(6)$.

For a different approach to this classification see [BKL2].

\section{Basic Definitions}

Let $A$ be an arbitrary (not necessarily associative) algebra over
$\mathbf{C}$.  By a formal distribution
$$a(z) = \sum_{i \in \mathbf{Z}} a(i) z^{-i-1} \in A[[z]]$$
we mean a power series over $A$, which is infinite in both
directions.

Two formal distributions $a(z)$, $b(z)$ are said to be mutually
local if there exists an integer $N = N(a,b) \geq 0$ such that
$a(z)b(w)(z-w)^N = b(w)a(z) (z-w)^N = 0$.

We will consider a countable family of operations:
$$a(z)\circ_n b(z) = Res_w a(w)b(z)(w-z)^n, \; n \geq 0, \; n
\in \mathbf{Z}.$$ Here $Res_w$ means the coefficient at $w^{-1}$.

If $a(z)$, $b(z)$ are mutually local then only finitely many
products $a \circ_n b$ may be different from zero.

\begin{defin}  A vector space $C \subseteq A[[z^{-1},z]]$ is
called a conformal algebra of formal distributions if $\partial C
\subseteq C$, $\partial = \frac{d}{dz}$, $C \circ_n C \subseteq C$
for an arbitrary $n \geq 0$ and every two elements from $C$ are
mutually local.
\end{defin}

By Dong Lemma (see [K2]) if $A$ is an associative or Lie algebra
then for an arbitrary collection $C$ of pairwise mutually local
distributions the closure of $C$ with respect to the action of
$\partial$ and to all operations $\circ_n$ , $n \geq 0$, is a
conformal algebra of formal distributions.

\begin{exam} (1) Let $\cal G$ be an arbitrary algebra and let $A =
{\cal G}[t, t^{-1}]$ be the algebra of Laurent polynomials over
$\cal G$.  For an arbitrary element $a \in {\cal G}$ let $\tilde a
= \sum_{i \in \mathbf{Z}} (at^i)z^{-i-1} \in A[[z^{-1},z]]$.

Any two formal distributions $\tilde a$, $\tilde b$ are mutually
local.

\m

(2) Let $\cal V$ir =$Der \mathbf{C}[t^{-1}, t]$ be the
(centerless) Virasoro algebra. The formal distribution
$$L = \sum_{i \in \mathbf{Z}} t^{i+1} \frac{d}{dt} z^{-i-2} \in {\cal
V}ir[[z^{-1}, z]]$$ is mutually local with itself.

\m

(3) Let $W = <t^{-1}, t, \frac{d}{dt}>$ be the (associative) Weyl
algebra of differential operators on $C[t^{-1},t]$.  Let $J_k =
\sum_{i \in Z} t^i ({\frac{d}{dt}})^k z^{-i-1}$, $k \geq 0$. Any
two formal distributions $J_k$, $J_l$ are mutually local.

\end{exam}

\m

In all three cases (1), (2) and (3) we can talk about the
conformal algebras $Cur({\cal G})$, ${\cal V}ir$, $W$
respectively, generated by them.

\m

Now we are ready to introduce an abstract definition of a
conformal algebra.

\m

Let $C$ be a module over a polynomial algebra
$\mathbf{C}[\partial]$, which is equipped with countably many
binary bilinear operations $C \circ_n C \rightarrow C$, $n \geq
0$.

\begin{defin}  We say that $(C, \partial, \circ_n)$ is an abstract
conformal algebra if for arbitrary elements $a,b \in C$ arbitrary
$n \geq 0$, we have:

1)  $\partial(a \circ_n b)= \partial a \circ_n b + a \circ_n
\partial b$,

2) $\partial a \circ_n b = - na \circ_{n-1} b$; for $n = 0$ the
condition turns into $\partial a \circ_0 b = 0$.

3) (Locality) There exists an integer $N = N(a,b) \geq 0$ such
that for an arbitrary $n \geq N$ we have $a \circ_n b = 0$.
\end{defin}

Every conformal algebra of formal distributions is an abstract
conformal algebra.  The converse is also true: every conformal
algebra can be realized as an algebra of formal distributions over
some algebra of coefficients.  Moreover, among these algebras of
coefficients there is a universal one $Coef\!f(C)$.

\begin{defin} We say that a conformal algebra $C$ is a Lie (resp.
associative, Jordan) algebra iff $Coef\!f(C)$ is a Lie (associative,
Jordan) algebra.
\end{defin}

Now let $C$ be a Lie conformal algebra and let $M$ be another
$\mathbf{C}[\partial]$-module.  Suppose that we have a family of
bilinear maps $C \circ_n M \subseteq M$, $n \geq 0$.

\begin{defin}  We say that $M$ is a conformal $C$-module if the
null split extension $C + M$ is a Lie conformal algebra.
\end{defin}

As above, $M$ can be realized as a space of formal distributions
over $Coef\!f(M)$, where $Coef\!f(M)$ is a universal (with this
property) Lie module over $Coef\!f(C)$.

\m

{\bf \underline{Important Remark} }  If there is a natural (and
standard) way to arrange elements of a (super)algebra $L$ in
formal distributions then we will talk about $L$ and modules over
$L$ even if we have in mind their conformal counterparts.

\section{The Cheng-Kac Superalgebra}

The exceptional conformal superalgebra $CK_6$ was introduced in
[CK2] and in [GLS].  In [MZ1] we constructed, for an arbitrary
associative commutative superalgebra $R$ with an even derivation
$d : R \rightarrow R$, a superalgebra $CK(R, d)$ so that $CK_6
\simeq CK(\mathbf{C}[t^{-1},t], \frac{d}{dt})$.

\m

Lets recall the construction of $CK(R,d)$ from [MZ1].

\m

Consider the associative Weyl algebra $W = \sum_{i \geq 0} Rd^i$,
where the variable $d$ does not commute with a coefficient $a \in
R$, but $da = ad + d(a)$.  We will realize the $CK(R,d)$ as a
superalgebra of $8 \times 8$ matrices over $W$.

\m

The simple finite dimensional Lie superalgebra $P(n-1)$ is the
superalgebra of $2n \times 2n$ matrices of the type $\left(
\begin{array}{cc}
a & k \\
h & -a^t
\end{array} \right),$ where $a,h,k$ are $n \times n$-matrices over
$\mathbf{C}$, $tr(a) = 0$, $k^t = -k$, $h^t = h$.  The
superalgebras $P(n)$, $n \neq 3$, are centrally closed.  However,
$P(3)$ has a nontrivial central cover $\hat{P(3)}$. Its existence
comes from the fact that the Lie algebra $K_4(\mathbf{C})$ of
skew-symmetric $4 \times 4$ matrices is a direct sum of two ideals
$K_4({\mathbf C}) = sl_2({\mathbf C}) \oplus sl_2({\mathbf C})$.
For an arbitrary element $k \in K_4(\mathbf{C})$ we consider its
decomposition $ k = k' + k''$ and let $\varphi(k) = k' - k''$. The
universal central cover $\hat{P(3)}$ of $P(3)$ can be realized as
a superalgebra of $8 \times 8$-matrices over the polynomial
algebra ${\mathbf C}[d]$ of the type
$$\left(\begin{array}{cc}
a & k \\
\varphi(k)d + h & -a^t
\end{array} \right) + \alpha d I_8,$$
where $a, k, h$ are $4 \times 4$ matrices over $\mathbf{C}$,
$tr(a) =0$, $k = -k^t$, $h = h^t$, $\alpha \in \mathbf{C}$ and
$I_8$ is the identity matrix.

\m

The superalgebra $CK(R,d)$ is a subsuperalgebra of $8 \times 8$
matrices over $W$ generated by $\hat{P(3)}$ and by all matrices
$\left(
\begin{array}{cc}
e_{ij}(a) & 0 \\
0 & -e_{ji}(a)
\end{array} \right)$ where $a \in R$, $1 \leq i \neq j
\leq 4$.

\m

The Cartan subalgebra $H$ of $CK(R,d)$ consists of diagonal
matrices
$$H = \{h = diag(a_1, \ldots, a_4, -a_1, \ldots, -a_4), \;
a_i \in {\mathbf C}, \; \sum_{i=1}^4 a_i = 0 \},$$ the even and
the odd roots of the $CK(R,d)$ with respect to the action of $H$
are:
$$\Delta_{\bar0} = \{w_i - w_j \, | \, 1 \leq i \neq j \leq 4 \},$$
$$ \Delta_{\bar1} = \{ w_i + w_j, \; 1 \leq i \neq j \leq 4, \;
-w_i-w_j, \; 1 \leq i,j \leq 4 \}.$$
Notice that $w_i(a) = a_i$, $
1 \leq i \leq 4$.

Thus, the superalgebra $CK(R,d)$ is graded by the abelian group
$$\sum_{i=1}^4 {\mathbf{Z}}w_i/{\mathbf{Z}}(w_1 + w_2+w_3 +
w_4),$$
 $CK(R,d) = \sum_{\alpha \in \Delta \cup \{0\}}
CK(R,d)_{\alpha}$.

\m

Let us fix the notation for the following weight elements:

$$ e _{w_i - w_j} = \left(
\begin{array}{cc}
e_{ij}& 0 \\
0 & -e_{ji}
\end{array} \right), \; e_{w_i - w_j}(a) = \left(
\begin{array}{cc}
e_{ij}(a) & 0 \\
0 & -e_{ji}(a)
\end{array} \right),$$

$$ h _{w_i - w_j}(a) = \left(
\begin{array}{cc}
e_{ii}(a) - e_{jj}(a) & 0 \\
0 & e_{jj}(a) - e_{ii}(a)
\end{array} \right), \; q_{-w_i - w_j} = \left(
\begin{array}{cc}
0 & 0 \\
e_{ij} + e_{ji} & 0
\end{array} \right),$$

$$ q _{-w_i - w_j}(a) = \left(
\begin{array}{cc}
0 & 0 \\
 e_{ij}(a) + e_{ji}(a)& 0
\end{array} \right), \; q_{w_i + w_j} = \left(
\begin{array}{cc}
0 & e_{ij} - e_{ji} \\
\varphi(e_{ij} - e_{ji})d & 0
\end{array} \right), $$

$a \in R$.

\bigskip

In [MZ3] it was shown that $CK(R,d)_{w_i - w_j} = e_{w_i
-w_j}(R)$, $1 \leq i \neq j \leq 4$;  $CK(R,d)_{-2w_i} =
q_{-2w_i}(R)$; $CK(R,d)_{w_i + w_j} = [q_{w_i + w_k}, e_{w_j -
w_k}(R)] + q_{-w_k - w_l}(R)$, where $\{i,j,k,l\} = \{1,2,3,4\}$.

For an arbitrary element $a \in R$ consider the element
$$[[e_{w_4 -w_1}(a), q_{w_3 + w_1}], q_{w_2 + w_1}] =$$

$$ \left(
\begin{array}{cc}
e_{11}(da) + e_{22}(ad) + e_{33}(ad) + e_{44}(ad) & 0 \\
 0  & e_{11}(ad) + e_{22}(da) + e_{33}(da) + e_{44}(da)
\end{array} \right) $$

$$= I_8(ad) - \left(
\begin{array}{cc}
e_{11}(a') & 0 \\
 0 & -e_{11}(a') + I_4(a')
\end{array} \right), \; a' = [a,d] = d(a)$$

We will denote the element on the right hand side as $Vir(a)$. The
mapping  $ad \rightarrow Vir(a)$ from $Rd \rightarrow {\cal
V}ir(R)$ is and isomorphism of Lie algebras.

\m

It was shown in [MZ3] that $CK(R,d)_0 = H \otimes R + {\cal
V}ir(R)$.

Consider the functional
$$f : \sum_{i = 1}^4 {\mathbf{Z}}w_i / {\mathbf{Z}}(\sum_{i=1}^4 w_i)
\rightarrow \mathbf{Z}$$ given by  $f(w_1) = 5, \, f(w_2) = -3, \,
f(w_3) = 2, \, f(w_4) = -4$.

Notice that $f(\pm w_i \pm w_j) \neq 0$, unless $\pm w_i \pm w_j =
0$.

From now on we will denote $L = CK_6 = CK({\mathbf{C}}[t^{-1},t],
\frac{d}{dt})$.

\m

Note that $L_0 = H \otimes {\mathbf{C}}[t^{-1},t] > \! \! \! \lhd
{\cal V}ir(R) \leq Cur(sl_4) > \! \! \! \lhd {\cal V}ir(R)  \leq
L$.

\m

The algebra $L$ has a triangular decomposition $L = L_{-} + L_0 +
L_{+}$, \newline $L_{-} = \sum_{f(\alpha) < 0} L_{\alpha}$, $L_{+}
= \sum_{f(\alpha) > 0} L_{\alpha}$.

\m

Let $M$ be a conformal module of finite type over the Lie
conformal algebra $CK_6$.  Then the subalgebra $sl_4 \subseteq L$
acts on $M$ and the action of $sl_4$ commutes with the action of
the polynomial algebra ${\mathbf{C}}[\partial]$. Hence $M$
decomposes into a finite direct sum of eigenspaces with respect to
the action of $H$, $$M = \sum_{\gamma \in H^*} M_{\gamma}.$$

If $M$ is irreducible, then there exists a unique highest weight
$\lambda \in H^*$ such that $M_{\lambda} \neq (0)$ and $L_{+}
\circ_n M_{\lambda} = (0)$ for all $n \geq 0$; $M_{\lambda}$ is an
irreducible conformal module over $L_0$.

We have mentioned above that $L_0 \subset Cur(sl_4) + {\cal V}ir
\subset L$.

\m

Let $M'$ be the $Cur(sl_4) + {\cal V}ir$-module generated by
$M_{\lambda}$.  Let $M''$ be the largest submodule of $M'$ such
that $M'' \cap M_{\lambda} = 0$.  Then $M'/M''$ is an irreducible
$Cur(sl_4) + {\cal V}ir$-module and $(M'/M'')_{\lambda} =
M_{\lambda}$.  Let $V = Coef\!f(M)$ be an $L$-module.

From the description of irreducible modules of finite type over
\break $Cur (sl_4) > \! \! \! \lhd {\cal V}ir(R)$  (see  [CK1]) it
follows that the module $V_{\lambda}$ can be identified with
${\mathbf{C}}[t^{-1},t]$, say $V_{\lambda} =
\overline{{\mathbf{C}}[t^{-1},t]}$.  For arbitrary elements $a,b
\in {\mathbf{C}}[t^{-1},t]$, $h \in H$ we have $(h \otimes a)\bar
b = <\lambda, h> \overline{ab}$.  Moreover, there exist scalars
$\alpha, \beta \in \mathbf{C}$ such that for arbitrary $a,b \in
{\mathbf{C}}[t^{-1},t]$ we have $$Vir(a) \bar b = \overline{-ab' +
\beta a'b + \alpha ab}.$$

\m

Denote this $L_0$-module as $V(\lambda, \beta, \alpha)$. It is
well known that, for an arbitrary $\lambda \in H^*$,
given an irreducible $L_0$-module $W$ such that the elements $h
\in H$ act on $W$ as scalar multiplications $<\lambda, h>$, there
exists a unique $L$-module with the highest weight $\lambda$ under
the action of $H$, whose  $\lambda$-space is isomorphic to $W$ as
$L_0$-module.  If we consider the irreducible $L_0$-module
$V(\lambda, \beta, \alpha)$, then the corresponding irreducible
$L$-module will be denoted as $Irr(\lambda, \beta, \alpha)$.

\m

It follows from the above that every irreducible conformal module
over $CK_6$ is isomorphic to $Irr(\lambda, \beta, \alpha)$ for
some $\lambda \in H^*$, $\beta, \alpha \in \mathbf{C}$.  This
gives rise to the question:

\b

{\it For which parameters $\lambda \in H^*$, $\beta, \alpha \in
\mathbf{C}$, the irreducible conformal module $Irr(\lambda,
\Delta, \alpha)$ is of finite type?}

\b

Let $\lambda$ be an integral dominant weight, that is, $<\lambda,
w_1 - w_3>$, $<\lambda, w_3 - w_2 >$, $<\lambda, w_2 - w_4 >$ all
lie in ${\mathbf{Z}}_{\geq 0}$.

\begin{teor}  The conformal module $Irr(\lambda, \beta, \alpha)$
is of finite type if and only if

(1)$<\lambda, h_{w_1 - w_3}> \geq 2$; $\beta, \alpha \in
\mathbf{C}$, or

(2)$<\lambda, h_{w_1 - w_3}> =1$; $<\lambda, h_{w_2-w_3}> = 0$,
$\beta = -1$, $\alpha \in \mathbf{C}$.

These modules exhaust all conformal irreducible $CK_6$-modules of
finite type.

\end{teor}

Since $V(\lambda, \beta, \alpha)$ are known to be conformal
modules of finite type over $L_0$ (see [CK1]), we can easily
conclude that $Irr(\lambda, \beta, \alpha)$  is of finite type if
and only if it has finitely many weights with respect to the
action of $H$. At this point we can forget about conformal modules
and address the question:

\b

{\it For which $\lambda \in H^*$, $\beta, \alpha \in \mathbf{C}$, the
$L$-module $Irr(\lambda, \beta, \alpha)$ has finitely many
weights?}

\begin{lema} Let $\alpha = w_i - w_j$ or $-w_i - w_j$.  For an
element $a \in R$, let $X_{\alpha}(a) = e_{w_i-w_j}(a)$ or
$q_{-w_i - w_j}(a)$ defined as above.  Suppose that $\alpha < 0$
and for any decomposition $- \alpha = \alpha_1 + \cdots +
\alpha_r$ into a sum of positive roots, for any elements $x_i \in
L_{\alpha_i}$, $1 \leq i \leq r$,there exist an element $h \in H$
and an element $b \in R$ such that  $[x_1, [x_2, \ldots [x_r,
X_{\alpha}(a)] \cdots] = h \otimes ab$ for an arbitrary $a \in R$.
Then for an arbitrary element $v_{\lambda} \in V_{\lambda}$ we
have $X_{\alpha}(a)v_{\lambda} = X_{\alpha}(1)(av)_{\lambda}$.
\end{lema}

{\bf Proof:} It is sufficient to show that
$$U(L_+)(X_{\alpha}(a)v_{\lambda} -
X_{\alpha}(1)(av)_{\lambda})\cap V_{\lambda} = (0).$$

Otherwise there exists a decomposition $- \alpha = \alpha_1 +
\cdots + \alpha_r$, $\alpha_i > 0$ and elements $x_i \in
L_{\alpha_i}$ such that  $x_1 \cdots x_r(X_{\alpha}(a)v_{\alpha} -
X_{\alpha}(1)(av)_{\lambda}) \neq 0$.

But  $$x_1 \cdots x_r X_{\alpha}(a)v_{\lambda} = [x_1, [x_2,
\ldots , [x_r, X_{\alpha}(a)] \ldots]v_{\lambda} = (h \otimes
ab)v_{\lambda} =$$
$$ h(abv)_{\lambda} = [x_1, [x_2, \ldots , [x_r, X_{\alpha}(1)] \ldots](av)_{\lambda}
= x_1 \ldots x_r X_{\alpha}(1)(av)_{\lambda},$$ a contradiction.
The lemma is proved.

\begin{lema}  The negative roots $w_2 - w_3$, $w_4 - w_3$, $-w_1 -
w_2$, $-w_1 - w_3$, $- w_1 - w_4$, $w_4 - w_2$ satisfy the
assumptions of Lemma 3.1
\end{lema}

{\bf Proof:}  We list all possible decompositions.  The roots $w_3
- w_2$, $w_1 + w_4$ and $w_2 - w_4$ do not have nontrivial
decompositions.  Then, $w_1 + w_2 = (w_1 + w_4) + (w_2 - w_4)$,
$w_1 + w_3 = (w_1 + w_2) + (w_3 - w_2) = (w_1 + w_4) + (w_2 - w_4)
+ (w_3 - w_2) = (w_3 - w_4) + (w_4 + w_1)$; $w_3 - w_4 = (w_2 -
w_4)+ (w_3 - w_2)$.

\m

The condition of Lemma 3.1 is checked by a straightforward
computation in the superalgebra $L$.  The lemma is proved.

\begin{lema}  For arbitrary elements $a,b \in R$ we have
$$ [q_{w_1+w_4}, e_{w_2 - w_4}(a)][q_{w_1 + w_4}, e_{w_3 -
w_4}(b)]q_{-w_1 - w_2} q_{-w_1 - w_3} v_{\lambda} = \gamma
(abv)_{\lambda},$$ where $\gamma = <\lambda, h_{w_1 - w_2}>(1 -
<\lambda, h_{w_1 - w_3}>)$.
\end{lema}

{\bf Proof:}  We have  $[q_{w_1 + w_4}, e_{w_3 - w_4}(b)]q_{-w_1 -
w_2}q_{-w_1-w_3}v_{\lambda} = (I) - (II)$, where
$$(I) = q_{w_1 + w_4} e_{w_3 - w_4}(b)q_{-w_1 - w_2} q_{-w_1 - w_3}
v_{\lambda}=$$
$$ q_{w_1 + w_4}q_{-w_1 - w_2} e_{w_3 - w_4}(b) q_{-w_1 - w_3}
v_{\lambda} =$$
$$  q_{w_1 + w_4}q_{-w_1 - w_2} [e_{w_3 - w_4}(b), q_{-w_1 - w_3}]
 v_{\lambda} +
  q_{w_1 + w_4}q_{-w_1 - w_2} q_{-w_1 - w_3} e_{w_3 - w_4}(b)
v_{\lambda} =$$
$$ -q_{w_1 + w_4}q_{-w_1 - w_2} q_{-w_1 - w_4}(b)v_{\lambda} =
-q_{w_1 + w_4}q_{-w_1 - w_2} q_{-w_1 - w_4}(bv)_{\lambda},$$ by
Lemma 3.1.

$$(II) = e_{w_3 - w_4}(b) q_{w_1 + w_4}q_{-w_1 - w_2}q_{-w_1 -
w_3}v_{\lambda} =$$
$$ - e_{w_3 - w_4}(b) e_{w_4 - w_2}q_{-w_1 - w_3}v_{\lambda}- e_{w_3 - w_4}(b)q_{-w_1 -
w_2} q_{w_1 + w_4} q_{-w_1 - w_3}v_{\lambda} =$$
$$ - e_{w_3 - w_4}(b) q_{-w_1 - w_3}e_{w_4 - w_2}v_{\lambda}- q_{-w_1 -w_2} e_{w_3 - w_4}(b)q_{w_1 + w_4} q_{-w_1 - w_3}v_{\lambda}.$$

Now  $$e_{w_3 - w_4}(b)q_{-w_1 - w_3}e_{w_4 - w_2}v_{\lambda} =$$
$$[e_{w_3 - w_4}(b),q_{-w_1 - w_3}]e_{w_4 - w_2}v_{\lambda} + q_{-w_1 - w_3}e_{w_3 - w_4}(b)e_{w_4 -
w_2}v_{\lambda}.$$

The second summand is $0$ since $f(w_3 - w_4) = 6$ whereas $f(w_4
- w_2) = -1$.

The first summand is
$$ - q_{-w_1 - w_4}(b)e_{w_4 - w_2}v_{\lambda} = $$
$$ - [q_{-w_1 - w_4}e_{w_4 - w_2}]v_{\lambda} - e_{w_4 -w_2}q_{-w_1 -
w_4}(b)v_{\lambda}=$$
$$ - q_{-w_1 - w_2}(b)v_{\lambda} - e_{w_4 -w_2}q_{-w_1 -
w_4}(b)v_{\lambda}=$$
$$ - q_{-w_1 - w_2}(bv)_{\lambda} - e_{w_4 -w_2}q_{-w_1 -
w_4}(bv_{\lambda}),$$ by Lemma 3.2.

\m

As for the other summand of (II),

$$ q_{-w_1 -w_2} e_{w_3 - w_4}(b)q_{w_1 + w_4} q_{-w_1 -
w_3}v_{\lambda}= $$
$$q_{-w_1 -w_2} [e_{w_3 - w_4}(b),[q_{w_1 + w_4}, q_{-w_1 -
w_3}]]v_{\lambda}= $$
$$q_{-w_1 -w_2} h_{w_3 - w_4}(b)v_{\lambda} = q_{-w_1 -w_2} h_{w_3 -
w_4}(bv)_{\lambda}.$$

We proved that $$[q_{w_1 + w_4}, e_{w_3-w_4}(b)]q_{-w_1
-w_2}q_{-w_1 -w_3}v_{\lambda} = P(bv)_{\lambda},$$  where $P$ is
an operator that does not involve $b$.  Choosing $b = 1$ we get
$$P = ad([q_{w_1 + w_4}, e_{w_3-w_4}]) ad(q_{-w_1-w_2}) ad(q_{-w_1
-w_3}) =$$ $$ ad(q_{w_3+w_1}) ad(q_{-w_1-w_2}) ad(q_{-w_1
-w_3}).$$

Now we have to consider the element
$$[q_{w_1+w_4}, e_{w_2 - w_4}(a)]q_{w_3 + w_1}q_{-w_1 -
w_2}q_{-w_1-w_3}(bv)_{\lambda}.$$

Remark that $[[q_{w_1+w_4}, e_{w_2 - w_4}(a)], q_{w_3 + w_1}] \in
e_{w_1 - w_4}(R)$ and  $$e_{w_1 - w_4}(R)q_{-w_1 -
w_2}q_{-w_1-w_3}(bv)_{\lambda} = [[ e_{w_1 - w_4}(R),q_{-w_1 -
w_2}], q_{-w_1-w_4}](bv)_{\lambda} = (0).$$

Hence our expression becomes

$$ - q_{w_1+w_3} [q_{w_1+ w_4}, e_{w_2 - w_4}(a)]q_{-w_1 -w_2}q_{-w_1-w_3}(bv)_{\lambda}.$$

Denote  $X =[q_{w_1+ w_4}, e_{w_2 - w_4}(a)]$, $Y = q_{-w_1
-w_2}$, $Z = q_{-w_1-w_3}$.  Then $XYZ = [[X,Y],Z] - Y[X,Z] + YZX
+ Z[X,Y]$, $[X,Y] = h_{w_2 - w_1}(a)$, $[X,Z] = e_{w_2 - w_3}(a)$,
$[[X,Y],Z] = q_{-w_1 - w_3}(a)$.  By Lemma 2, $$h_{w_2 -
w_1}(a)(bv)_{\lambda} = h_{w_2 - w_1}(abv)_{\lambda}, \; e_{w_2 -
w_3}(a)(bv)_{\lambda} = e_{w_2 - w_3}(abv)_{\lambda}$$ and
$$q_{-w_1-w_3}(a)(bv)_{\lambda} = q_{-w_2 -
w_3}(abv)_{\lambda}.$$

As we did above, we can conclude that
$$[q_{w_1 + w_4}, e_{w_2 - w_4}(a)][q_{w_1+w_4}, e_{w_3 -
w_4}(b)]q_{-w_1-w_2}q_{-w_1 - w_3}v_{\lambda} = \tilde
P(abv)_{\lambda},$$ where $\tilde P$ is an operator that does not
involve $a$ or $b$.  Choosing $a = b = 1$ we get
$$\tilde P = ad(q_{w_2 + w_1})ad(q_{w_3 + w_1})ad(q_{-w_1 -
w_2})ad(q_{-w_1 - w_3}).$$

Now $$\tilde P (abv)_{\lambda} = q_{w_2 + w_1}[q_{w_3 + w_1},
q_{-w_1 - w_2}]q_{-w_1 - w_3}(abv)_{\lambda} - $$ $$q_{w_2 +
w_1}q_{-w_1-w_2}q_{w_3+w_1}q_{-w_1 -w_3}(abv)_{\lambda}=$$
$$[q_{w_2 + w_1}, [[q_{w_3 + w_1}, q_{-w_1
-w_2}],q_{-w_1-w_3}]](abv)_{\lambda} - $$
$$[q_{w_2 + w_1},  q_{-w_1-w_2}],[ q_{w_3 + w_1},
q_{-w_1-w_3}]](abv)_{\lambda}= $$
$$- h_{w_2-w_1}(abv)_{\lambda} - h_{w_2 - w_1}h_{w_3 -
w_1}(abv)_{\lambda}= \gamma (abv)_{\lambda},$$  $\gamma =
<\lambda, h_{w_1-w_2}>(1-<\lambda, h_{w_1-w_3}>)$.  The lemma is
proved.

\b

{\bf Remark.}  In what follows  $[x_1, \ldots, x_n]$ denotes the
left-normed commutator $[...[x_1,x_2],x_3], \ldots, x_n]$.

\begin{lema} $[[[e_{w_4 - w_1}(a), [q_{w_1+w_4},e_{w_3 -
w_4}(c)]],[q_{w_1 + w_4}, e_{w_2 - w_4}(b)]] =  \newline h_{w_1 -
w_4}(ab'c) - Vir(abc)$.
\end{lema}

{\bf Proof:} \, Since $[e_{w_4 - w_1}(a), q_{w_1+w_4}] = 0$ the
left hand side is equal to $$[e_{w_4 - w_1}(a), e_{w_3 - w_4}(c),
q_{w_1+w_4}, e_{w_2 - w_4}(b),q_{w_1+w_4} ] =$$ $$ -[ e_{w_3 -
w_1}(ac), q_{w_1+w_4},  e_{w_2 - w_4}(b), q_{w_1+w_4}].$$

Furthermore, $e_{w_2 - w_4}(b) = -[q_{w_1+w_2}, q_{-w_1-w_4}(b)]$.
Substituting this expression we get:

$$LHS = [e_{w_3 - w_1}(ac), q_{w_1+w_4},
q_{w_1+w_2},q_{-w_1-w_4}(b), q_{w_1+w_4}]+$$
$$[e_{w_3-w_1}(ac), q_{w_1+w_4}, q_{-w_1-w_4}(b), q_{w_1+w_2}, q_{w_1+w_4}] = (I) + (II).$$

Recall that we use the notation $[e_{w_4 - w_1}(a), q_{w_3 + w_1},
q_{w_2+w_1}] = Vir(a)$. Now,
$$[e_{w_3-w_1}(ac), q_{w_1+w_4},q_{w_1 + w_2}] =
-[e_{w_4-w_1}(ac),e_{w_3 -w_4}, q_{w_1 + w_4},q_{w_1+w_2}] =$$ $$
-[e_{w_4-w_1}(ac),[e_{w_3-w_4}, q_{w_1 + w_4}], q_{w_1+w_2}]$$
since $[e_{w_4 -w_1}(ac),q_{w_1+w_4}] \subseteq L_{2w_4} = (0)$.

Using that $[e_{w_3-w_4},q_{w_1+w_4}] = -q_{w_3+w_1}$, our
expression becomes
$$ [e_{w_4-w_1}(ac), q_{w_3+w_1},q_{w_2+w_1}] = Vir(ac).$$

Hence,
$$(I) = [e_{w_3-w_1}(ac), q_{w_1+w_4},
q_{w_1+w_2},[q_{-w_1-w_4}(b), q_{w_1+w_4}]] =$$
$$-[Vir(ac),h_{w_1-w_4}(b)] = -h_{w_1-w_4}(ab'c).$$

On the other hand,
$$(II) = [e_{w_3 - w_1}(ac), q_{w_1+w_4},
q_{-w_1-w_4}(b),q_{w_1+w_2},q_{w_1+w_4}] =$$
$$[e_{w_3-w_1}(ac),[q_{w_1+w_4},q_{-w_1-w_4}(b)],q_{w_1+w_2},q_{w_1+w_4}]
= $$
$$[e_{w_3-w_1}(ac),h_{w_1-w_4}(b),q_{w_1+w_2},q_{w_1+w_4}] =
[e_{w_3-w_1}(abc),q_{w_1+w_2},q_{w_1+w_4}]=$$
$$- [e_{w_3-w_1}(abc),q_{w_1+w_4},q_{w_1+w_2}] = -Vir(abc),$$
as we have seen above.  This proves the lemma.

\m

Lemma 3.4 implies that
$$[q_{w_1+w_4},
e_{w_2-w_4}(b)][q_{w_1+w_4},e_{w_3-w_4}(c)]e_{w_4-w_1}(a)v_{\lambda}
= $$
$$ - [e_{w_4-w_1}(a),[q_{w_1+w_4},e_{w_3-w_4}(c)],[q_{w_1+w_4},
e_{w_2-w_4}(b)]]v_{\lambda} =$$
$$(-h_{w_1-w_4}(ab'c)-Vir(abc))v_{\lambda} =$$
$$ - <\lambda,h_{w_1-w_4}> (abcv')_{\lambda} +((\mu_0 abc + \mu_1a'bc + \mu_2ab'c
+\mu_3abc')v)_{\lambda},$$ here $v$ is viewed as an element from
$R = F[t^{-1},t]$; $\mu_0, \; \mu_1, \; \mu_2, \; \mu_3$ are
scalars from $F$.

Choosing $a = 1$ we get:
$$[q_{w_1+w_4},e_{w_2-w_4}(b)][q_{w_1+w_4},e_{w_3-w_4}(c)]
e_{w_4-w_1}v_{\lambda} = $$  $$((\mu_0 bc +  \mu_2b'c
+\mu_3bc')v)_{\lambda} - <\lambda,h_{w_1-w_4}> (bcv')_{\lambda}.$$

Hence, $[q_{w_1+w_4},e_{w_2-w_4}(b)][q_{w_1+w_4},e_{w_3-w_4}(c)]
e_{w_4-w_1}(av)_{\lambda}=$
$$((\mu_0 abc +  \mu_2 ab'c
+\mu_3 abc')v)_{\lambda} - <\lambda,h_{w_1-w_4}>
((abcv')_{\lambda} + (a'bcv)_{\lambda}).$$

This implies

$$[q_{w_1+w_4},e_{w_2-w_4}(b)][q_{w_1+w_4},e_{w_3-w_4}(c)](
e_{w_4-w_1}(a)(v)_{\lambda} - e_{w_4-w_1}(av)_{\lambda} =$$
$$(\mu_1 + <\lambda,h_{w_1-w_4}>) (a'bcv)_{\lambda}.  \hskip 3 true cm  (*)$$

\section{ The case $ <\lambda, h_{w_1-w_3}> \: \geq 2$}

In this section we will prove that if $\lambda$ is an integral
dominant functional and $<\lambda, h_{w_1 - w_3}>  \geq 2$, then
for arbitrary $\beta, \alpha \in F$ the irreducible module
$V(\lambda, \beta, \alpha)$ has only finitely many weights with
respect to the action of $H$.

\medskip

Let $\gamma = <\lambda, h_{w_1 - w_2}>(1 - <\lambda, h_{w_1 -
w_3}>)$ (see Lemma 3.3).  Since $w_1 - w_2 = (w_1 - w_3)+(w_3 -
w_2)$ and the root $w_3 - w_2$ is positive, we conclude that
$<\lambda,  h_{w_1 - w_2}> \geq 2$ and therefore $\gamma \neq 0$.
Let $\xi = \frac{\mu_1 + <\lambda,h_{w_1-w_4}>}{\gamma}$.

\begin{lema}  $e_{w_4 - w_1}(a)v_{\lambda} - e_{w_4 -
w_1}(av)_{\lambda} - \xi q_{-w_1 - w_2}q_{-w_1 -
w_3}(a'v)_{\lambda} = 0$.

\end{lema}

{\bf Proof.} Denote the left hand side of the above equality as
$w$. In order to prove that $w = 0$ we need only to check that
$L_{+}w \cap V_{\lambda} = (0)$.  From the equality (*) and Lemma
3.3 it follows that $$[q_{w_1 + w_4}, e_{w_2 - w_4}(b)][q_{w_1 +
w_4}, e_{w_3 - w_4}(c)]w = 0.$$

Consider the element $\; q_{-w_3 - w_4}(b)w$.  We have  $\;$
$q_{-w_3 -w_4}(b)e_{w_4 - w_1}(a)v_{\lambda} =$\newline $q_{-w_1 -
w_3}(ab)v_{\lambda} = q_{-w_1 - w_3}(abv)_{\lambda}$ by Lemma 3.2.
The last expression is equal also to $q_{-w_3 - w_4}(b)e_{w_3 -
w_1}(av)_{\lambda}$. Furthermore, $$q_{-w_3 - w_4}(b)q_{-w_1
-w_2}q_{-w_1 - w_3}v_{\lambda} = q_{-w_1 - w_2}q_{-w_1 -
w_3}q_{-w_3 - w_4}(b) v_{\lambda} = 0,$$ since $- w_3 - w_4$ is
positive.

\medskip

We have shown that $q_{-w_3 - w_4}(R)w = (0)$.

Similarly  $q_{-w_2 - w_4}(R)w = (0)$.

\medskip

Let us  show that $e_{w_1 - w_4}(R)w = (0)$.  Indeed, $f(w_1 -
w_4) = 9$, $f(-w_1 - w_2) = -2, f(-w_1 - w_3) = -7$.  Hence,
$$e_{w_1 - w_4}(b)q_{-w_1 - w_2}q_{-w_1 - w_3}v_{\lambda} = [e_{w_1
- w_4}(b), q_{-w_1 - w_2}, q_{-w_1 - w_3}]v_{\lambda} = 0.$$  Now
$$e_{w_1 - w_4}(b)w = e_{w_1 - w_4}(b)e_{w_4 - w_1}(a)v_{\lambda}
-e_{w_1 - w_4}(b)e_{w_4 - w_1}(av)_{\lambda} =$$ $$ h_{w_1 -
w_4}(ab)v_{\lambda} - h_{w_1 - w_4}(b)(av)_{\lambda} = 0.$$

\medskip

Since $[L_{w_1+w_2}, L_{w_1 + w_3}] \subseteq e_{w_1 - w_4}(R)$ it
follows that for arbitrary elements  $x \in L_{w_1 + w_2}$, $y \in
L_{w_1 + w_3}$, $xyw = -yxw$.

\medskip

We have  $L_{w_1+w_2} = [q_{w_1 + w_4}, e_{w_2 - w_4}(R)] +
q_{-w_3 - w_4}(R)$, $L_{w_1+w_3} = [q_{w_1 + w_4}, e_{w_3 -
w_4}(R)] + q_{-w_2 - w_4}(R)$.

\medskip

From what we proved above, it follows that $L_{w_1+w_2}L_{w_1 +
w_3}w = (0)$.  Together with $e_{w_1 - w_4}(R)w = (0)$ it implies
that $U(L_+)w \cap V_{\lambda} = (0)$ and therefore $w = 0$. Lemma
is proved.

\begin{lema}  $e_{w_3 - w_1}(a)v_{\lambda} = e_{w_3 -
w_1}(av)_{\lambda} - \xi q_{-w_1 -w_2}q_{-w_1 -
w_4}(a'v)_{\lambda}$.
\end{lema}

{\bf Proof.}  By Lemma 2.5, $$e_{w_3 - w_4}e_{w_4 -
w_1}(a)v_{\lambda} = e_{w_3 - w_4}e_{w_4 - w_1}(av)_{\lambda} +
\xi e_{w_3 - w_4}q_{-w_1 - w_2}q_{-w_1 - w_3}(a'v)_{\lambda}.$$

Since $w_3 - w_4$ is positive, it implies
$$[e_{w_3 - w_4}, e_{w_4 -w_1}(a)]v_{\lambda} = [e_{w_3 - w_4}, e_{w_4
-w_1}](av)_{\lambda}+ $$  $$\xi q_{-w_1 - w_2}[e_{w_3 - w_4},
q_{-w_1 - w_3}](a'v)_{\lambda} = e_{w_3 - w_1}(av)_{\lambda} - \xi
q_{-w_1 - w_2}q_{-w_1 - w_4}(a'v)_{\lambda}.$$ This implies the
result and proves the lemma.

From now on in this section, unless otherwise stated, we will
assume that $<\lambda, h_{w_1 - w_3}> = 2$.   Our first aim is to
show that $e^3_{w_3 - w_1} V_{\lambda} = (0).$

\begin{lema}
$e_{w_1 - w_3}(a)e^3_{w_3 - w_1} v_{\lambda} = 6 \, \xi \, e_{w_3
- w_1}q_{-w_1 - w_2}q_{-w_1 - w_4}(a'v)$.

\end{lema}

{\bf Proof.}  Taking into account that $$[e_{w_1 - w_3}(a), e_{w_3
- w_1}, e_{w_3 - w_1}, e_{w_3 -w_1}] = e_{w_1 - w_3}(a)e^3_{w_3 -
w_1} - $$  $$3 e_{w_3 - w_1} e_{w_1 - w_3}(a)e^2_{w_3 - w_1} + 3
e^2_{w_3 - w_1} e_{w_1 - w_3}(a)e_{w_3 - w_1} - e^3_{w_3 - w_1}
e_{w_1 - w_3}(a) = 0$$ and  $$e_{w_1 - w_3}(a)e^2_{w_3 - w_1} =
[e_{w_1 - w_3}(a), e_{w_3 - w_1}, e_{w_3 - w_1}] + $$ $$ 2e_{w_3 -
w_1}e_{w_1 - w_3}(a) e_{w_3 - w_1} - e^2_{w_3 - w_1}e_{w_1 -
w_3}(a)=$$ $$ -2 e_{w_3 - w_1}(a) + 2 e_{w_3 - w_1} e_{w_1 -
w_3}(a) e_{w_3 - w_1} - e^2_{w_3 - w_1} e_{w_1 - w_3}(a),$$ we get
$$e_{w_1 - w_3}(a) e^3_{w_3 - w_1} v_{\lambda} = 3e_{w_3 - w_1} (-2 e_{w_3 -
w_1}(a)+ 2 e_{w_3 - w_1} e_{w_1 - w_3}(a)e_{w_3 - w_1})v_{\lambda}
-$$ $$ 3 e^2_{w_3 - w_1} e_{w_1 - w_3}(a) e_{w_3 - w_1}
v_{\lambda} = -6 e_{w_3 - w_1}e_{w_3 - w_1}(a) v_{\lambda} + 3
e^2_{w_3 - w_1}h_{w_1 - w_3}(a)v_{\lambda}= $$ $$ -6e_{w_3 - w_1}
(e_{w_3 - w_1}(av)_{\lambda} - \xi q_{-w_1-w_2} q_{-w_1 -
w_4}(a'v)_{\lambda} + 6 e^2_{w_3 - w_1}(av)_{\lambda}=$$ $$6 \,
\xi e_{w_3 - w_1} q_{- w_1 - w_2}q_{-w_1 - w_4}(a'v)_{\lambda}.$$
The lemma is proved.

\begin{lema}  $e_{w_1 - w_3}(a) e_{w_1-w_3}(b)e^3_{w_3 -
w_1}v_{\lambda} = 0$.
\end{lema}

{\bf  Proof.}  By Lemma 4.3, the left hand side is equal to $$6 \,
\xi e_{w_1 - w_3}(a) e_{w_3 - w_1} q_{-w_1-w_2} q_{-w_1 -
w_4}(b'v)_{\lambda}.$$

We notice that $$e_{w_1 - w_3}(a)q_{-w_1 -
w_2}q_{-w_1-w_4}(b'v)_{\lambda} = [e_{w_1 -
w_3}(a),q_{-w_1-w_2},q_{-w_1-w_4}](b'v)_{\lambda} = 0.$$

Hence, \hskip 2 true cm  $e_{w_1-w_3}(a) e_{w_3 - w_1} q_{-w_1 -
w_2}q_{-w_1 - w_4}(b'v)_{\lambda} = $  $$[e_{w_1 - w_3}(a), e_{w_3
- w_1}]q_{-w_1 - w_2}q_{-w_1 - w_4}(b'v)_{\lambda} =  h_{w_1 -
w_3}(a) q_{-w_1 - w_2} q_{-w_1 - w_4}(b'v)_{\lambda}.$$

From Lemma 3.2 it follows that the last expression is equal to
$$h_{w_1 - w_3}(1)q_{-w_1 - w_2}q_{-w_1 - w_4}(ab'v)_{\lambda} =$$
$$<\lambda+ w_3 - w_1, h_{w_1 - w_3}(1)> q_{-w_1 - w_2}q_{-w_1 -
w_4}(ab'v)_{\lambda} = 0.$$ The lemma is proved.

\begin{lema} $ [q_{w_1 + w_3}, e_{w_4 - w_3}(a)]e_{w_3 -
w_1}q_{-w_1 - w_2}q_{-w_1 - w_4}v_{\lambda} =$

$q_{w_4 + w_1}e_{w_3 - w_1} q_{-w_1 - w_2} q_{-w_1 -
w_4}(av)_{\lambda}$.
\end{lema}

{\bf Proof.}  We have  $[q_{w_1 + w_3}, e_{w_4 - w_3}(a)] = q_{w_1
+ w_3} e_{w_4 - w_3}(a) - e_{w_4 - w_3}(a)q_{w_1 + w_3}$. Now,
since the total weight of the expression  $ q_{w_1 + w_3} e_{w_3 -
w_1}q_{-w_1 - w_2}q_{-w_1 - w_4}$ is $3w_3 - w_1$ that is
positive, we only need to consider the expression
$$q_{w_1+w_3}e_{w_4 - w_3}(a)e_{w_3 - w_1} q_{-w_1 - w_2} q_{-w_1
- w_4} v_{\lambda} = (I) + (II)$$ where $(I) = q_{w_1+w_3} e_{w_4
- w_1}(a)q_{-w_1-w_2}q_{-w_1-w_4}v_{\lambda}$,  and \newline $(II)
= q_{w_1 + w_3}e_{w_3 - w_1}e_{w_4 - w_3}(a)q_{-w_1-w_2}q_{-w_1 -
w_4} v_{\lambda}$.

\medskip

Let us consider these expressions separately.

$$(I) = q_{w_1 + w_3}q_{-w_1-w_2} e_{w_4 - w_1}(a) q_{-w_1 - w_4}
v_{\lambda} = $$  $$\underbrace{q_{w_1+w_3} q_{-w_1 -
w_2}q_{-w_1-w_4}e_{w_4-w_1}(a)v_{\lambda}}_{I.1} +
\underbrace{q_{w_1 + w_3}q_{-w_1 - w_2}[e_{w_4 - w_1}(a), q_{-w_1
- w_4}]v_{\lambda}}_{I.2}.$$

$$(I.1) = q_{w_1 + w_3}q_{-w_1 - w_2}q_{-w_1 - w_4} e_{w_4 -
w_1}(a)v_{\lambda} = $$
$$q_{w_1 + w_3}q_{-w_1 - w_2}q_{-w_1 - w_4}(e_{w_4 -
w_1}(av)_{\lambda} + \xi q_{-w_1 - w_2} q_{-w_1 -
w_3}(a'v)_{\lambda}) =$$
$$q_{w_1 + w_3}q_{-w_1 - w_2}q_{-w_1 - w_4}e_{w_4 -
w_1}(av)_{\lambda},$$ because  $q_{-w_1-w_2}q_{-w_1 - w_4}q_{-w_1
- w_2} = 0$;

$$(I.2) = - q_{w_1+w_3}q_{-w_1 - w_2}q_{-2w_1}(a)v_{\lambda} =$$
$$q_{-w_1-w_2}q_{w_1+w_3}q_{-2w_1}(a)v_{\lambda} - [q_{w_1 +
w_3}, q_{-w_1 - w_2}]q_{-2w_1}(a)v_{\lambda} = (I.2.1) +
(I.2.2);$$

$$(I.2.1) = q_{-w_1 - w_2}[q_{w_1+w_3}, q_{-2w_1}(a)]v_{\lambda} =
-2q_{-w_1 -w_2}e_{w_3-w_1}(a) v_{\lambda} =$$
$$ -2q_{-w_1 - w_2}(e_{w_3 - w_1}(av)_{\lambda} - \xi q_{-w_1 -
w_2}q_{-w_1 - w_4}(a'v)_{\lambda}) = -2 q_{-w_1 - w_2}e_{w_3 -
w_1}(av)_{\lambda};$$

$$(I.2.2) =  e_{w_3 - w_2}q_{-2w_1}(a)v_{\lambda} =
q_{-2w_1}(a)e_{w_3 - w_2}v_{\lambda} = 0$$ since $w_3 - w_2$ is
positive.

$$(II) = q_{w_1+w_3}e_{w_3 - w_1}e_{w_4 - w_3}(a)q_{-w_1 -
w_2}q_{-w_1 - w_4}v_{\lambda} =$$
$$q_{w_1 + w_3}e_{w_3 - w_1}q_{-w_1 - w_2}e_{w_4 -
w_3}(a)q_{-w_1-w_4}v_{\lambda} =$$
$$\underbrace{q_{w_1 + w_3}e_{w_3 - w_1}q_{-w_1 - w_2}q_{-w_1 -
w_4}e_{w_4 - w_3}(a)v_{\lambda}}_{II.1} +$$
$$\underbrace{q_{w_1 + w_3}e_{w_3 - w_1}q_{-w_1 - w_2}[e_{w_4 -
w_3}(a), q_{-w_1 - w_4}]v_{\lambda}}_{II.2}.$$

But  $$(II.1) = q_{w_1 + w_3}e_{w_3 - w_1}q_{-w_1 - w_2}q_{-w_1 -
w_4}e_{w_4 - w_3}(a)v_{\lambda} =$$
$$q_{w_1 + w_3}e_{w_3 - w_1}q_{-w_1 - w_2}q_{-w_1
- w_4}e_{w_4 - w_3}(av)_{\lambda}$$ by Lemma 3.2;

$$(II.2) = -q_{w_1 + w_3}e_{w_3 -
w_1}q_{-w_1-w_2}q_{-w_1-w_3}(a)v_{\lambda} =$$ $$
-q_{w_1+w_3}e_{w_3 - w_1}q_{-w_1 -
w_2}q_{-w_1-w_3}(av)_{\lambda}$$ again by Lemma 3.2.

To summarize, we have proved that  $$[q_{w_1 + w_3}, e_{w_4 -
w_3}(a)]e_{w_3 - w_1}q_{-w_1 - w_3}q_{-w_1 - w_4}v_{\lambda} =
P(av)_{\lambda},$$ where $P$ is an operator that does not involve
$a$.  Choosing $a = 1$, we get $P = ad(q_{w_4 + w_1})ad(e_{w_3 -
w_1})ad(q_{-w_1 - w_2})ad(q_{-w_1-w_4})$.  The lemma is proved.

\begin{lema}  (i)  $[q_{w_1 + w_3},e_{w_4 - w_3}(a)]e^3_{w_3 -
w_1}v_{\lambda} =$  $$ -q_{w_1+w_4}e^3_{w_3-w_1}(av)_{\lambda} -3
\,\xi e^2_{w_3-w_1}q_{-w_1 - w_2}(a'v)_{\lambda}.$$

(ii) $[q_{w_1 + w_3},e_{w_4 - w_3}(a)][q_{w_1 + w_3},e_{w_4 -
w_3}(b)]e^3_{w_3 - w_1}v_{\lambda} =$ $$3 \,\xi q_{w_1+w_4}
e^2_{w_3 - w_1} q_{-w_1 - w_2}((ab' - a'b)v)_{\lambda}.$$

(iii) $[q_{w_1 + w_3},e_{w_4 - w_3}(a)][q_{w_1 + w_3},e_{w_4 -
w_3}(b)][q_{w_1 + w_3},e_{w_4 - w_3}(c)]e^3_{w_3 - w_1}v_{\lambda}
=0.$

\end{lema}

{\bf Proof.} (i)  The element $q_{w_1 + w_3}$ commutes with
$e_{w_3 - w_1}$ and $w_1 + w_3$ is positive.  Hence,
$q_{w_1+w_3}e^3_{w_3 - w_1}v_{\lambda} = 0.$  Furthermore,
$[e_{w_4 - w_3}(a), e_{w_3 - w_1},e_{w_3 - w_1}]= 0$.  Hence,
$$e_{w_4 - w_3}(a) e^3_{w_3 - w_1}v_{\lambda} = 3 \, e^2_{w_3 - w_1}e_{w_4 - w_3}(a)e_{w_3 -
w_1}v_{\lambda}- 2 \, e^3_{w_3 - w_1} e_{w_4 - w_3}(a) v_{\lambda}
= $$ $$3 \, e^2_{w_3 - w_1}e_{w_4 - w_1}(a)v_{\lambda} + e^3_{w_3
- w_1} e_{w_4 - w_3}(a) v_{\lambda} =$$  $$3 \, e^2_{w_3 - w_1}(
e_{w_4 - w_1}(av)_{\lambda} + \xi q_{-w_1 - w_2}q_{-w_1 -
w_3}(a'v)_{\lambda}) + e^3_{w_3 - w_1} e_{w_4 -
w_3}(av)_{\lambda}.$$

We have proved that  $[q_{w_1+w_3}, e_{w_4 - w_3}(a)]e^3_{w_3 -
w_1}v_{\lambda} =$  $$[q_{w_1+w_3}, e_{w_4 - w_3}]e^3_{w_3 -
w_1}(av)_{\lambda} + 3 \,\xi q_{w_1+w_3}e^2_{w_3 - w_1}
q_{-w_1-w_2}q_{-w_1 - w_3}(a'v)_{\lambda}=$$
$$-q_{w_1+w_4}e^3_{w_3 - w_1}(av)_{\lambda} + 3 \,\xi e^2_{w_3 - w_1}q_{w_1+w_3}
q_{-w_1-w_2}q_{-w_1 - w_3}(a'v)_{\lambda} =$$
$$ -q_{w_1+w_4}e^3_{w_3 - w_1}(av)_{\lambda} -3 \,\xi e^2_{w_3 - w_1}e_{w_3 - w_2}q_{-w_1 - w_3}(a'v)_{\lambda} -$$
$$3 \,\xi e^2_{w_3 - w_1}q_{-w_1-w_2}q_{w_1+w_3}q_{-w_1 -
w_3}(a'v)_{\lambda}=$$
$$-q_{w_1+w_4}e^3_{w_3 - w_1}(av)_{\lambda} + 3 \,\xi e^2_{w_3 - w_1}q_{-w_1 -
w_2}(a'v)_{\lambda} - 3 \,\xi e^2_{w_3 - w_1}q_{-w_1 -
w_2}h_{w_1-w_3}(a'v)_{\lambda} = $$
$$-q_{w_1+w_4}e^3_{w_3 - w_1}(av)_{\lambda} - 3 \,\xi e^2_{w_3 - w_1}q_{-w_1 -
w_2}(a'v)_{\lambda}.$$ The assertion (i) is proved.

\medskip

(ii) Let us apply (i) to $[q_{w_1 + w_3},e_{w_4 - w_3}(b)]e^3_{w_3
- w_1}v_{\lambda}$  and consider both summands of the right hand
side of (i) separately.

\medskip

We have, $[q_{w_1 + w_3},e_{w_4 - w_3}(a)]q_{w_1+w_4}e^3_{w_3 -
w_1}(bv)_{\lambda}= $
$$-q_{w_1 + w_4}[q_{w_1 + w_3},e_{w_4 -
w_3}(a)]e^3_{w_3 - w_1}(bv)_{\lambda} =$$
$$q_{w_1+w_4}(q_{w_1+w_4}e^3_{w_3-w_1}(abv)_{\lambda} + 3 \, \xi e^2_{w_3-w_1}q_{-w_1 -
w_2}(a'bv)_{\lambda}) =$$
$$ 3 \, \xi q_{w_1+w_4} e^2_{w_3-w_1}q_{-w_1 -
w_2}(a'bv)_{\lambda}).$$

Acting on the second summand, we get
$$[q_{w_1+w_3}, e_{w_4 - w_3}(a)]e^2_{w_3-w_1}q_{-w_1 -
w_2}(b'v)_{\lambda} =$$
$$ -e^2_{w_3-w_1}[q_{w_1+w_3}, e_{w_4 - w_3}(a)]q_{-w_1 -
w_2}(b'v)_{\lambda})\; + $$  $$2 e_{w_3-w_1}[q_{w_1+w_3}, e_{w_4 -
w_3}(a)]e_{w_3-w_1}q_{-w_1 - w_2}(b'v)_{\lambda} =$$
$$e^2_{w_3-w_1}[q_{w_1+w_3}, e_{w_4 -
w_3}(a)]q_{-w_1 - w_2}(b'v)_{\lambda} \; + $$ $$2
e_{w_3-w_1}[q_{w_1+w_3}, e_{w_4 - w_3}(a), e_{w_3-w_1}]q_{-w_1 -
w_2}(b'v)_{\lambda} = $$
$$e^2_{w_3-w_1}[q_{w_1+w_3},
e_{w_4 - w_3}(a), q_{-w_1-w_2}](b'v)_{\lambda} \; +$$
$$ 2 e_{w_3-w_1}[q_{w_1+w_3},
e_{w_4 - w_1}(a)]q_{-w_1 - w_2}(b'v)_{\lambda} =$$
$$ e^2_{w_3-w_1}e_{w_4 - w_2}(a)(b'v)_{\lambda} \, +  \, 2e_{w_3-w_1}[q_{w_1+w_3},
e_{w_4 - w_1}(a)]q_{-w_1 - w_2}(b'v)_{\lambda}.$$

The first summand of this sum is equal to $ e^2_{w_3-w_1}e_{w_4 -
w_2}(ab'v)_{\lambda}$ by Lemma 3.2.  As for the second summand,
$$ e_{w_3-w_1}[q_{w_1+w_3}, e_{w_4 -
w_1}(a)]q_{-w_1-w_2}(b'v)_{\lambda} = $$
$$ e_{w_3-w_1}q_{w_1+w_3}e_{w_4 -
w_1}(a)q_{-w_1-w_2}(b'v)_{\lambda}= $$
$$ e_{w_3-w_1}q_{w_1+w_3}q_{-w_1-w_2}e_{w_4 -
w_1}(a)(b'v)_{\lambda} =$$
$$ e_{w_3-w_1}q_{w_1+w_3}q_{-w_1-w_2}e_{w_4-w_1}(ab'v)_{\lambda} +$$
$$ \xi
e_{w_3-w_1}q_{w_1+w_3}\underbrace{q_{-w_1-w_2}q_{-w_1-w_2}}_0q_{-w_1-w_3}(a'b'v)_{\lambda}=
$$
$$e_{w_3-w_1}q_{w_1+w_3}q_{-w_1-w_2}e_{w_4-w_1}(ab'v)_{\lambda}.
$$

We have shown that
$$[q_{w_1+w_3},
e_{w_4-w_3}(a)]e^2_{w_3-w_1}q_{-w_1-w_2}(b'v)_{\lambda} =
P(ab'v)_{\lambda},$$ where $P$ is an operator which does not
involve $a$ or $b$.  Choosing $a = 1$, $b = t$, we get  $P =
-ad(q_{w_1+w_4})ad(e_{w_3-w_1})^2 ad(q_{-w_1-w_2})$ which finishes
the proof of (ii).

\bigskip

(iii)  By using (ii) we need only to show that
$$[q_{w_1+w_3},e_{w_4-w_3}(a)]q_{w_1+w_4}e^2_{w_3-w_1}q_{-w_1-w_2}v_{\lambda}
= 0.$$  Since  $[q_{w_1+w_3}, e_{w_4-w_3}(a)]$  and $q_{w_1+w_4}$
commute, the expression above is
\newline $-q_{w_1+w_4}[q_{w_1+w_3},e_{w_4 -
w_3}(a)]e^2_{w_3-w_1}q_{-w_1-w_2}v_{\lambda}.$

We proved above that $$[q_{w_1+w_3},e_{w_4 -
w_3}(a)]e^2_{w_3-w_1}q_{-w_1-w_2}v_{\lambda} = -
q_{w_1+w_4}e^2_{w_3-w_1}q_{-w_1-w_2}(av)_{\lambda}.$$ Now
multiplying this expression on the left by $q_{w_1+w_4}$ we get 0.
This concludes the proof of the lemma.

\begin{lema} $[q_{w_1+w_3}, e_{w_2 -
w_3}(R)]^3e^3_{w_3-w_1}v_{\lambda} = (0)$.

\end{lema}

{\bf Proof.}  Apply  $ad{(e_{w_2-w_4})}^3$ to the equality
$[q_{w_1+w_3}, e_{w_4 - w_3}(R)]^3e^3_{w_3-w_1}v_{\lambda} = (0)$
of Lemma 4.6(iii).

Since $[e_{w_2-w_4}, q_{w_1+w_3}] = [e_{w_2-w_4}, e_{w_3-w_1}] =
[e_{w_4-w_3}(R), e_{w_2-w_4}, e_{w_2-w_4}] = (0)$, we will get
$$[q_{w_1+w_3}, [e_{w_2-w_4}, e_{w_4 -
w_3}(R)]]^3e^3_{w_3-w_1}v_{\lambda} = (0),$$ completing the proof
of the lemma.

\begin{lema} $e^3_{w_3-w_1}v_{\lambda} = 0$.

\end{lema}

{\bf Proof.}  If $e^3_{w_3-w_1}v_{\lambda} \neq 0$, then there
exist positive roots  $\alpha_1, \ldots, \alpha_s$ such that  $(0)
\neq L_{\alpha_1} \cdots L_{\alpha_s}e^3_{w_3-w_1}v_{\lambda}
\subseteq V_{\lambda}$.  Let $s$ be the minimal number with this
property.  Since we can move each $L_{\alpha_i}$ to the right
modulo shorter products, we can assume that for each $i$, $1 \leq
i \leq s$, $\alpha_i + w_3 - w_1$ is a root or 0 and $f(\alpha_i)
\leq f(w_1 - w_3) = 3$.  Among all positive roots, only $w_1 -
w_3$, $w_1+w_2$, $w_1+w_4$ have this properties.  Suppose that
$$(0) \neq L^i_{w_1+w_2}L^j_{w_1+w_4}e_{w_1-w_3}(R)^k e^3_{w_3-w_1}
v_{\lambda} \subseteq V_{\lambda}.$$   Then $i(w_1+w_2) +
j(w_1+w_4) + (3-k)(w_3-w_1) = m(w_1+w_2+w_3+w_4)$; $0 \leq i,j,k
\in \mathbf Z$, $m \in \mathbf{Z}$.

\medskip

This implies  $i = j = m$, $k = 3-m \geq 0$.  Hence, we have 3
options:

1)  $k = 2$ or $3$.  This contradicts Lemma 4.4.

\medskip

2)  $k = 1$.  By Lemma 4.3
$$L^2_{w_1+w_2}L^2_{w_1+w_4}e_{w_1-w_3}(R)e^3_{w_3-w_1}v_{\lambda}
\subseteq
L^2_{w_1+w_2}L^2_{w_1+w_4}e_{w_3-w_1}q_{-w_1-w_2}q_{-w_1-w_4}v_{\lambda}.$$

The factors in $L^2_{w_1+w_2}L^2_{w_1+w_4}$ on the right hand side
anticommute, because of the minimality of $s$ and the fact that
$[L_{w_1+w_2}, L_{w_1+w_4}] \subseteq e_{w_1 - w_3}(R)$, which
leads to the case 2).

\medskip

Suppose that at least one of the two $L_{w_1+w_4}$ factors lies in
$q_{-w_2 - w_3}(R)$.  Then
$$\underbrace{q_{-w_2-w_3}(a)e_{w_3-w_1}}q_{-w_1-w_2}q_{-w_1-w_4}v_{\lambda}
=$$
$$e_{w_3-w_1}q_{-w_2-w_3}(a)q_{-w_1-w_2}q_{-w_1-w_4}v_{\lambda} +
q_{-w_2-w_1}(a)q_{-w_1-w_2}q_{-w_1-w_4}v_{\lambda}.$$

The first summand is 0 because $-w_2-w_3$ is positive.  The second
summand is equal to
$$q_{-w_1-w_2}q_{-w_1-w_4}q_{-w_2-w_1}(a)v_{\lambda} =
q_{-w_1-w_2}q_{-w_1-w_4}q_{-w_2-w_1}(av)_{\lambda}$$ by Lemma 3.2.
Now it remains to notice that
$q_{-w_1-w_2}q_{-w_1-w_4}q_{-w_1-w_2}= 0$.

Thus, we can assume that both factors from $L_{w_1+w_4}$ are
$[q_{w_1+w_3}, e_{w_4-w_3}(a_i)]$, $i = 1,2$.

\medskip

By Lemma 4.5 we have
$$[q_{w_1+w_3}, e_{w_4-w_3}(a_1)][q_{w_1+w_3},
e_{w_4-w_3}(a_2)]e_{w_3-w_1}q_{-w_1-w_2}q_{-w_1-w_4}v_{\lambda} =
$$
$$[q_{w_1+w_3},
e_{w_4-w_3}(a_1)]q_{w_4+w_1}e_{w_3-w_1}q_{-w_1-w_2}q_{-w_1-w_4}(a_2v)_{\lambda}.
$$

The element $[q_{w_1+w_3}, e_{w_4-w_3}(a_1)]$ anticommutes with
$q_{w_4+w_1}$.  Hence again by Lemma 4.5  $$[q_{w_1+w_3},
e_{w_4-w_3}(a_1)]q_{w_4+w_1}e_{w_3
-w_1}q_{-w_1-w_2}q_{-w_1-w_4}(a_2v)_{\lambda} =$$
$$-q_{w_4+w_1} q_{w_4+w_1} e_{w_3
-w_1}q_{-w_1-w_2}q_{-w_1-w_4}(a_1a_2v)_{\lambda} = 0.$$

\medskip

3) $k = 0$.  We have to examine
$L^3_{w_1+w_2}L^3_{w_1+w_4}e^3_{w_3-w_1}v_{\lambda}.$  As above,
we conclude that factors from $L_{w_1+w_2}$ and $L_{w_1+w_4}$
anticommute module the previous cases ($k \geq 1$).

\medskip

From $q_{-w_3-w_4}(R)^2e^3_{w_3-w_1}v_{\lambda} = (0)$, it follows
that no more than one factor from $L_{w_1+w_4}$ lies in
$q_{-w_2-w_3}(R)$.

On the other hand, Lemma 4.6 (iii) implies that exactly one factor
from $L_{w_1 +w_4}$ lies in $q_{-w_2-w_3}(R)$. Similarly,
$q_{-w_3-w_4}(R)^2e^3_{w_3-w_1}v_{\lambda} = (0)$ and Lemma 4.7
imply that exactly one factor from $L_{w_1+w_2}$ lies in
$q_{-w_3-w_4}(R)$.

Now we need to show that
$$[q_{w_1+w_3},
e_{w_2-w_3}(a_1)][q_{w_1+w_3},e_{w_2-w_3}(a_2)][q_{w_1+w_3},
e_{w_4-w_3}(b_1)][q_{w_1+w_3}, e_{w_4-w_3}(b_2)]$$
$q_{-w_3-w_4}(c_1)q_{-w_2-w_3}(c_2)e^3_{w_3-w_1}v_{\lambda} = 0.$

\medskip

First notice that $q_{-w_2-w_3}(c_2)e^3_{w_3-w_1}v_{\lambda} =$
$$3e^2_{w_3-w_1}[q_{-w_2-w_3}(c_2),e_{w_3-w_1}]v_{\lambda} =$$
$$3e^2_{w_3-w_1}q_{-w_1-w_2}(c_2)v_{\lambda} = 3e^2_{w_3-w_1}q_{-w_1-w_2}(c_2
v)_{\lambda}$$  by Lemma 3.2.  Hence, without loss of generality,
we can assume that $c_2 = 1$ and similarly, $c_1 = 1$. Moreover,
$$q_{-w_3-w_4}q_{-w_2-w_3}e^3_{w_3-w_1}v_{\lambda} =$$
$$6e_{w_3-w_1}[q_{-w_3-w_4}, e_{w_3-w_1}][q_{-w_2-w_3},
e_{w_3-w_1}] v_{\lambda} =
6e_{w_3-w_1}q_{-w_1-w_4}q_{-w_1-w_2}v_{\lambda}.$$

Now,
$$[q_{w_1+w_3},e_{w_4-w_3}(b_1)][q_{w_1+w_3},e_{w_4-w_3}(b_2)]e_{w_3-w_1}q_{-w_1-w_4}q_{-w_1-w_2}v_{\lambda}
=0$$  follows from Lemma 4.5.  This concludes the proof of the
lemma.

\begin{lema} 1)  $(e_{w_4-w_3})^{<\lambda,h_{w_3-w_4}>+1}
v_{\lambda}= 0$,

\bigskip

2) $(e_{w_4-w_2})^{<\lambda,h_{w_2-w_4}>+1} v_{\lambda}= 0$.

\end{lema}

{\bf Proof.}  The only positive roots $\alpha$ such that $\alpha +
w_4 - w_3$ is a root and $f(\alpha) \leq f(w_3-w_4) = 6$ are $w_3
- w_4$, $w_1+w_2$, $w_3-w_2$, $w_2-w_4$.  Suppose that  $$(0) \neq
L^i_{w_1+w_2} L^j_{w_3-w_2}
L^k_{w_2-w_4}L^l_{w_3-w_4}e^q_{w_4-w_3} v_{\lambda} \subseteq
V_{\lambda}, \; \:  q = <\lambda, h_{w_3-w_4}>+1.$$  Then $
i(w_1+w_2) + j(w_3-w_2) + k(w_2-w_4)+(q-l)(w_4-w_3) =
m(w_1+w_2+w_3+w_4)$, where $i,j,k,l,m \in {\mathbf Z}$.

 This implies
$i = m, i-j+k = m, j-(q-l) = m$, $-k+(q-l) = m$.  Hence $i = m =
0$, $j = k = q-l \geq 0$.

\medskip

Now we have to examine the expression $ L^j_{w_3-w_2}
L^j_{w_2-w_4}L^l_{w_3-w_4}e^q_{w_4-w_3} v_{\lambda}$, where $l = q
- j$.

\medskip

Suppose that $l \geq 1$.  There exist rational numbers $\mu, \;
\nu$ such that for an arbitrary element $a \in R$,
$$e_{w_3-w_4}(a)e_{w_4-w_3}^q v_{\lambda} = \mu e^{q-2}_{w_4-w_3}
e_{w_4-w_3}(a)v_{\lambda} + \nu e^{q-1}_{w_4-w_3}h_{w_3-w_4}
v_{\lambda} = $$
$$\mu e^{q-2}_{w_4-w_3}e_{w_4-w_3}(av)_{\lambda} + \nu
e^{q-1}_{w_4-w_3}<\lambda, h_{w_3-w_4}>(av)_{\lambda}$$ by Lemma
3.2.

This implies that  $$e_{w_3-w_4}(a) e^{q}_{w_4-w_3}v_{\lambda} =
e_{w_3-w_4} e^{q}_{w_4-w_3}(av)_{\lambda} = 0$$ since $q =
<\lambda, h_{w_3-w_4}> + 1$.

\medskip

Now let $l = 0, j = q$.  As above
$$e_{w_2-w_4}(a)e^q_{w_4-w_3}v_{\lambda} = q e^{q-1}_{w_4-w_3}
e_{w_2-w_3}(a)v_{\lambda} = q
e^{q-1}_{w_4-w_3}e_{w_2-w_3}(av)_{\lambda}.$$  This implies that
$$L_{w_3-w_2}^q L_{w_2-w_4}^q e^q_{w_4-w_3} V_{\lambda} = L_{w_3-w_2}^q e_{w_2-w_4}^q e^q_{w_4-w_3}
V_{\lambda}$$ and, similarly, this expression is equal to
$e_{w_3-w_2}^q e_{w_2-w_4}^q e^q_{w_4-w_3} V_{\lambda}$.  We have
shown above that $$ e_{w_2-w_4}^q e^q_{w_4-w_3} v_{\lambda} = q!
e^q_{w_2-w_3} v_{\lambda}.$$ Now $e_{w_3-w_2}^q e_{w_2-w_3}^q
v_{\lambda} = 0$ because  $q = <\lambda, h_{w_3-w_4}>+1 \geq
<\lambda, h_{w_3-w_2}> + 1$.  This proves 1).  Let us prove now
assertion 2). The only positive roots $\alpha$ such that $\alpha
+w_4 - w_2$ is a root and $f(\alpha) \leq f(w_2 - w_4) = 1$ are
$w_2 - w_4$ and $w_1 + w_4$. If  $$ (0) \neq L_{w_1+w_4}^i
L_{w_2-w_4}^j e_{w_4-w_2}^p v_{\lambda} \subseteq V_{\lambda}, \:
p = <\lambda, h_{w_2-w_4}>+1,$$ then $i(w_1+w_4) + (p-j)(w_4-w_2)
= m(w_1+w_2+w_3+w_4)$, which implies $i = m = 0$, $j = p$. Arguing
as above, we see that $L_{w_2-w_4}^p e_{w_4-w_2}^p V_{\lambda} =
e_{w_2-w_4}^p e_{w_4-w_2}^p V_{\lambda} = 0$.  This completes the
proof of the lemma.

\begin{lema}  Let $M \subseteq L$ a subspace such that $M^n
v_{\lambda} = (0)$, where $v_{\lambda} \in V_{\lambda}$.  Let $1
\leq i \neq j \leq 4$ and $e^m_{w_i-w_j} v_{\lambda} = 0$. Suppose
further that $[e_{w_i-w_j}, M, M] = (0)$.  Then
$[M,e_{w_i-w_j}]^{m+n} v_{\lambda} = (0)$.

\end{lema}

{\bf Proof.}  From $[M,e_{w_i-w_j},e_{w_i-w_j}] = (0)$ it follows
that $$[M, e_{w_i-w_j}]^{m+n} = [ \underbrace{M \cdots M}_{m+n},
\underbrace{e_{w_i-w_j}, \ldots, e_{w_i-w_j}}_{m+n}],$$ where
products on the left hand side and in $M \cdots M$ are taken in
the associative algebra $End_F(V)$.  Hence, $$[M, e_{w_i -
w_j}]^{m+n} v_{\lambda} \subseteq \sum_{s+r = m+n} e_{w_i-w_j}^s
\underbrace{M \cdots M}e_{w_i - w_j}^r v_{\lambda}.$$ In each
nonzero summand on the right hand side  $\: r \leq m-1$.

From $[\underbrace{M, [M, [M, \ldots [M}_{r+1},e_{w_i-w_j}^r]]
\cdots ] = (0)$ it follows that  $$\underbrace{M \cdots M}_{m+n}
e_{w_i-w_j}^r \subseteq \sum_{p+q = m+n, p<m} M^p e_{w_i-w_j}^r
M^q$$ which implies that $q \geq n$ and therefore $M^q v_{\lambda}
= (0)$.

\begin{lema}  There exists $m \geq 1$ such that $e_{w_i-w_j}^m
V_{\lambda} = (0)$ for any $1 \leq i \neq j \leq m$.

\end{lema}

{\bf Proof.}  By Lemmas 4.8 and 4.9 the elements
$e_{\pm(w_1-w_3)}$, $e_{\pm(w_3-w_4)}$, $e_{\pm(w_2-w_4)}$ act
nilpotently on $V_{\lambda}$.  Now it remains to notice that those
elements generate $sl(4)$ and to use Lemma 4.10.

\begin{lema} For an arbitrary root $\alpha$ the subspace
$L_{\alpha}$ acts nilpotently on $V_{\lambda}$.

\end{lema}

{\bf Proof.}  Let $\alpha = w_i - w_j$, $1 \leq i \neq j \leq 4$,
$w_i - w_j$ negative.  We have shown that  $e_{w_i - w_j}^m
V_{\lambda} = (0)$.  Now, $L_{w_i-w_j} = [e_{w_j - w_i}(R), e_{w_i
- w_j},e_{w_i - w_j}]$  and $[e_{w_j - w_i}(R), e_{w_i -
w_j},e_{w_i - w_j}, e_{w_i - w_j}] = (0)$.  This implies that
$$L^m_{w_i-w_j} = [L^m_{w_j-w_i}, \underbrace{e_{w_i - w_j}, \ldots, e_{w_i -
w_j}}_{2m}] \subseteq \sum_{p+q=2m} e_{w_i - w_j}^p L^m_{w_j-w_i}
e_{w_i - w_j}^q.$$

If $q \geq m$ then $e^q_{w_i-w_j} V_{\lambda} = (0)$.  If $q \leq
m-1$, then $f(m(w_j - w_i) + q(w_i - w_j)) > 0$ and again
$L^m_{w_j - w_i} e^q_{w_i - w_j} V_{\lambda} = (0)$.  We have
shown that $L^m_{w_j - w_i}V_{\lambda} = (0)$.

\medskip

Let $\alpha$ be an odd root such that $L_{\alpha}$ acts on
$V_{\alpha}$ nilpotently, $\alpha$ is not of the form $-2w_k$.
Then for arbitrary $1 \leq i \neq j \leq 4$ the subspace
$[L_{\alpha}, e_{w_i-w_j}]$ acts on $V_{\lambda}$ nilpotently.
Indeed, since $\alpha \neq -2w_i$, we have $[L_{\alpha},
e_{w_i-w_j}, e_{w_i-w_j}] = [e_{w_i-e_j},L_{\alpha}, L_{\alpha}] =
(0)$.  Now the claim follows from Lemma 4.10.

Consider a root space $L_{w_i + w_j}$, $1 \leq i \neq j \leq 4$.
If one of $i,j$ is equal to 1, then $w_i + w_j > 0$.  Let $i \neq
1$, $j \neq 1$.  Then $[L_{w_i + w_j}, e_{w_i-w_1}] = (0)$, but
$<w_i+w_j, h_{w_i-w_1}> \neq 0$.  Hence, $L_{w_i+w_j} =
[[e_{w_i-w_1}, e_{w_1-w_i}],L_{w_i+w_j}] \subseteq [e_{w_i-w_1},
L_{w_1+w_j}].$

From what we proved above it follows that $L_{w_i + w_j}$ acts on
$V_{\lambda}$ nilpotently.

\medskip

Next, $L_{-2w_i} = [[e_{w_j-w_i}, e_{w_i-w_j}], L_{-2w_i}]
\subseteq [e_{w_j-w_i},L_{-w_i-w_j}]$, which implies that
$L_{-2w_i}$ acts on $V_{\lambda}$ nilpotently.  This completes the
proof of the lemma.

\section{Tensor product of modules $V(\lambda, \beta,\alpha)$}

In this section we will discuss a realization of modules $V(\beta,
\alpha)$ and define a tensor product in this class.

\medskip

Let $R$ be an arbitrary commutative $F$-algebra with a derivation
$d:R \rightarrow R$.  Recall that the Weyl algebra $W$ is $W =
\sum_{i = 0}^{\infty} Rd^i$, $da = ad + d(a)$.  For an arbitrary
scalar $\beta \in F$ consider the vector space $W_{\beta}(R,d) =
\{a_0d^{\beta} + a_1 d^{\beta - 1} + a_2 d^{\beta - 2} + \cdots ,
a_i \in R\}$, the (infinite) sums are understood formally, $\tilde
W(R,d) = \sum_{\beta \in F} W_{\beta}(R,d)$.  The rule $d^{\gamma}
a = \sum_{i=0}^{\infty} {\gamma \choose i} d^i(a)d^{\gamma - i}$,
where $d^i (a)$ is the $i$-th derivative of the element $a$, makes
$\tilde W(R,d)$ an associative algebra, $W \subseteq \tilde
W(R,d)$.  Moreover, for each $\beta \in F$ we have $[Rd,
W_{\beta}(R,d)] \subseteq W_{\beta}(R,d)$.  Hence $W_{\beta}(R,d)$
is a module over the Virasoro algebra $Rd$.

\medskip

Now consider the associative commutative algebra $\tilde R = R +
Rv$, $v^2 = 0$. Extend the derivation $d$ via $d(v) = -\alpha v$,
$\alpha \in F$.  Then the subspace $W_{\beta}(R, v, d) = \sum_{i =
0}^{\infty} Rvd^{\beta-i} \subset W_{\beta}(\tilde R, d)$ is an
$Rd$-submodule of $W_{\beta}(\tilde R, d)$.

\medskip

The following proposition is streightforward.

\begin{prop}  $W_{\beta}(R,v,d)/W_{\beta-1}(R,v,d) \simeq
V(\beta,\alpha)$.
\end{prop}

The tensor product $V(\beta_1,\alpha_1) \otimes_F
V(\beta_2,\alpha_2)$ can be identified with \break
$W_{\beta_1+\beta_2}(R,v_1v_2,d)/W_{\beta_1+\beta_2-1}(R,v_1v_2,d)$,
where $R = F[t_1^{-1},t_1, t_2^{-1}t_2]$, $d =$ \break $ -d/dt_1 -
d/dt_2$.

Since $d(t_1-t_2) = 0$ it follows that
$(t_1-t_2)(V(\beta_1,\alpha_1) \otimes V(\beta_2,\alpha_2))$ is a
submodule of $V(\beta_1,\alpha_1) \otimes V(\beta_2,\alpha_2)$.

Clearly, $V(\beta_1,\alpha_1) \otimes V(\beta_2,\alpha_2)/(t_1 -
t_2) \simeq V(\beta_1+\beta_2, \alpha_1+\alpha_2).$

\begin{prop}  If $V(\lambda_i,\beta_i,\alpha_i)$, $i = 1,2$ are
conformal modules of finite type, then so is  $ V(\lambda_1 +
\lambda_2, \beta_1+\beta_2,\alpha_1+\alpha_2)$.
\end{prop}

{\bf Proof.}  The $L$-modules $V(\lambda_i,\beta_i,\alpha_i)$ have
finitely many weight spaces with respect to the Cartan subalgebra
$H$ of $L$.  The tensor product $V = V(\lambda_1,\beta_1,\alpha_1)
\otimes V(\lambda_2,\beta_2,\alpha_2)$ also has finitely many
weight spaces.  The subspace of $V$ of weight
$\lambda_1+\lambda_2$ can be identified with $V(\beta_1,\alpha_1)
\otimes V(\beta_2,\alpha_2)$.Let $M$ be the submodule of $V$
generated by $(t_1 - t_2)(V(\beta_1,\alpha_1) \otimes
V(\beta_2,\alpha_2))$. Then $(V/M)_{\lambda_1+\lambda_2} \simeq
V(\beta_1 + \beta_2, \alpha_1+\alpha_2)$.  The $L$-module
$V(\lambda_1+\lambda_2,\beta_1 + \beta_2, \alpha_1+\alpha_2)$ is a
homomorphic image of the submodule of $V/M$ generated by
$(V/M)_{\lambda_1+\lambda_2}$.  Hence
$V(\lambda_1+\lambda_2,\beta_1 + \beta_2, \alpha_1+\alpha_2)$ has
finitely many weight spaces with respect to $H$.  This concludes
the proof of the proposition.

\vskip 1cm

Consider a copy of the algebra of Laurent polynomials
$\overline{F[t,t^{-1}]}$ and make it a $W$-module via $a \bar b =
\overline{ab}$, $d \bar b = - \bar b'$, $a,b \in F[t^{-1},t]$.
Then the space of 8-columns $\overline{F[t,t^{-1}]}^8$ becomes a
left module over $M_8(W)$, hence a $CK(6)$-module.  It is easy to
see that this $CK(6)$-module is irreducible.

\bigskip

If we define the form $(w_i/w_j) = \delta_{ij}$ on $\sum_{i=1}^4
Fw_i$ and view functionals on $H$ as elements of $\sum_{i=1}^4
Fw_i$, then the highest weight of the module
$\overline{F[t,t^{-1}]}^8$ is $w_1$, $(h_{w_i-w_j} \otimes a)(\bar
b, 0, \ldots, 0)^T = (w_i-w_j/w_1)(\bar b, 0, \ldots, 0)^T$.
Moreover $Vir(a) (\bar b, 0, \ldots, 0)^T = (\overline{-ab' -
a'b}, 0, \ldots, 0)^T$.  Hence $\overline{F[t,t^{-1}]}^8 \simeq
V(w_1,-1,0)$.

\begin{prop}  If $\lambda$ is an integral dominant functional and
\break $< \lambda, h_{w_1-w_3}> \geq 2$, then for arbitrary
$\beta, \alpha \in F$ the irreducible module \break $V(\lambda,
\beta, \alpha)$ has only finitely many weights with respect to the
action of $H$.

\end{prop}

{\bf Proof.}  Let $<\lambda, h_{w_1 - w_3}> = k \geq 2$. By
Proposition 5.2 the module $V' = V(\lambda - (k-2)w_1, \beta +
(k-2), \alpha)$ has finitely many $H$-weights.  Tensoring $V'$
with $\overline{F[t,t^{-1}]}^8 \simeq V(w_1,-1,0)$ $k-2$ times and
using Proposition 5.2 we get the result.

\section {The case $<\lambda, h_{w_1-w_3}> = 1$}

The aim of this section is to prove the following

\begin{prop} Let $\lambda$ be an integral dominant weight, such
that \break $ <\lambda, h_{w_1 - w_3}> = 1$.  Then $V(\lambda,
\beta, \alpha)$ has finitely many weights with respect to $H$ if
and only if $<\lambda, h_{w_3 - w_2}> = 0$ and $\beta = -1$.

\end{prop}

Suppose at first that $\lambda$ is an integral dominant weight
such that \break $<\lambda,h_{w_1-w_3}> = 1$ and $V(\lambda,
\beta, \alpha)$ has finitely many $H$-weights.

\begin {lema} For arbitrary elements $a \in R$, $v_{\lambda} \in
V_{\lambda}$ we have \break $e_{w_3-w_1}(a) v_{\lambda} =
e_{w_3-w_1}(av)_{\lambda}$.

\end{lema}

{\bf Proof.}  Since $V(\lambda, \beta, \alpha)$ is a finite sum of
eigenspaces with respect to the $H$ it follows that the element
$e_{w_3-w_1}$ acts on $V_{\lambda}$ nilpotently.  The standard
argument shows that $e^2_{w_3 - w_1}V_{\lambda} = (0)$.  Now for
an arbitrary $a \in R$ we have

$$ 0 = e_{w_1-w_3}(a)e_{w_3-w_1}^2v_{\lambda} =$$
$$[e_{w_1-w_3}(a), e_{w_3-w_1}, e_{w_3-w_1}]v_{\lambda} +
2e_{w_3-w_1} e_{w_1-w_3}(a)e_{w_3-w_1}v_{\lambda} =$$
$$  -2 e_{w_1-w_3}(a)v_{\lambda} + 2 e_{w_3
-w_1}h_{w_1-w_3}(a)v_{\lambda} =
 -2(e_{w_3-w_1}(a)v_{\lambda} - e_{w_1-w_3}(av)_{\lambda}).$$
This concludes the proof of the lemma.

\begin{lema}  $\beta = -1$.

\end{lema}

{\bf Proof.}  $$[q_{w_1+w_3},
e_{w_2-w_3}(b)][q_{w_1+w_3},e_{w_4-w_3}(c)]e_{w_3-w_1}(a)v_{\lambda}
=$$
$$[[q_{w_1+w_3},
e_{w_2-w_3}(b)],[[q_{w_1+w_3},e_{w_4-w_3}(c)],e_{w_3-w_1}(a)]]v_{\lambda}
=$$
$$(-h_{w_1-w_3}(ab'c) + Vir(abc))v_{\lambda}$$ as in Lemma 3.4.
This element is equal to $(-ab'cv - abcv' + \beta (abc)'v + \alpha
abcv)_{\lambda}$.

On the other hand, by Lemma 3.1 we have
 $$[q_{w_1+w_3},
e_{w_2-w_3}(b)][q_{w_1+w_3},e_{w_4-w_3}(c)]e_{w_3-w_1}(a)v_{\lambda}
=$$
$$[q_{w_1+w_3},
e_{w_2-w_3}(b)][q_{w_1+w_3},e_{w_4-w_3}(c)]e_{w_3-w_1}(av)_{\lambda}
=$$
$$(-h_{w_1-w_3}(b'c) + Vir(bc))(av)_{\lambda} = (-ab'cv -(bc)(av)' + \beta (bc)'(av) + \alpha abcv)_{\lambda}.$$

Comparing these two expressions we see that  $\beta a'bcv = -
bca'v$, so $\beta = -1.$  This finishes the proof of the lemma.

\begin{lema} $<\lambda, h_{w_3 - w_2}> = 0$.
\end{lema}

{\bf  Proof.} We have  $q_{w_1+w_2}q_{w_1 +
w_4}q_{-w_3-w_4}q_{-w_2-w_3}e^2_{w_3-w_1}v_{\lambda} = 0$.  Now,
$$q_{-w_2-w_3}e_{w_3-w_1}^2v_{\lambda} =
2e_{w_3-w_1}[q_{-w_2-w_3},e_{w_3-w_1}]v_{\lambda} =
2e_{w_3-w_1}q_{-w_1-w_2}v_{\lambda}.$$  Hence, \hskip 1.4 true cm
$0 =
q_{w_1+w_2}q_{w_1+w_4}q_{-w_3-w_4}e_{w_3-w_1}q_{-w_1-w_2}v_{\lambda}
=$
$$q_{w_1+w_4}q_{-w_3-w_4}q_{w_1+w_2}e_{w_3-w_1}q_{-w_1-w_2}v_{\lambda}
=$$
$$q_{w_1+w_4}q_{-w_3-w_4}[q_{w_1+w_2},e_{w_3-w_1}]q_{-w_1-w_2}v_{\lambda}
+$$
$$q_{w_1+w_4}q_{-w_3-w_4}e_{w_3-w_1}[q_{w_1+w_2},q_{-w_1-w_2}]v_{\lambda}=$$
$$q_{w_1+w_4}q_{-w_3-w_4}q_{w_2+w_3}q_{-w_1-w_2}v_{\lambda} +
<\lambda,h_{w_1-w_2}>q_{w_1+w_4}q_{-w_3-w_4}e_{w_3-w_1}v_{\lambda}=$$
$$q_{w_1+w_4}[q_{-w_3-w_4},q_{w_2+w_3}]q_{-w_1-w_2}v_{\lambda} -
q_{w_1+w_4}q_{w_2+w_3}q_{-w_3-w_4}q_{-w_1-w_2}v_{\lambda} +$$
$$<\lambda,h_{w_1-w_2}> [q_{w_1+w_4},[q_{-w_3-w_4},e_{w_3-w_1}]]v_{\lambda}
=$$
$$ - q_{w_1+w_4}q_{-w_1-w_4}v_{\lambda} +<\lambda,h_{w_1-w_2}> <\lambda,h_{w_1-w_4}> v_{\lambda} =$$
$$ <\lambda,h_{w_1-w_4}> ( -1 +
<\lambda,h_{w_1-w_2}>)v_{\lambda}.$$

Since  $<\lambda,h_{w_1-w_4}>  \geq <\lambda,h_{w_1-w_3}>  = 1$ it
follows that $<\lambda,h_{w_1-w_2}>  = 1$ and therefore
$<\lambda,h_{w_3-w_2}>  = 0$.  This concludes the proof of the
lemma.

\medskip

Now we will assume that $\lambda$ is an integral dominant weight
such that $<\lambda,h_{w_1-w_3}> = 1$, $<\lambda,h_{w_3-w_2}> =
0$.  Let $\beta = -1$.  We will prove that $V(\lambda, \beta,
\alpha) $ is a finite sum of eigenspaces with respect to $H$.

\begin{lema} Under the assumptions above, $e_{w_3 -
w_1}(a)v_{\lambda}= e_{w_3-w_1}(av)_{\lambda}$ for arbitrary $a
\in R$, $v_{\lambda} \in V_{\lambda}$.
\end{lema}

{\bf Proof.} The computations of Lemma 6.3 show that for
$<\lambda,h_{w_1-w_3}> = 1$, $\beta = -1$ we have $$
[q_{w_1+w_3},e_{w_2-w_3}(R)][q_{w_1+w_3},
e_{w_4-w_3}(R)](e_{w_3-w_1}(a)v_{\lambda} -
e_{w_3-w_1}(av))_{\lambda} = 0.$$

Also, $q_{-w_3-w_4}(R)(e_{w_3-w_1}(a)v_{\lambda}  -
e_{w_3-w_1}(av)_{\lambda}) = q_{-w_2-w_3}(R)(e_{w_3-w_1}(a)
v_{\lambda}  - e_{w_3-w_1}(av)_{\lambda}) = (0)$ by Lemma 3.2 .
This implies that  $U(L_+)(e_{w_3-w_1}(a)v_{\lambda} =
e_{w_3-w_1}(av)_{\lambda}$.  Lemma is proved.

\begin{lema} $e_{w_4 - w_1}(a)v_{\lambda} =
e_{w_4-w_1}(av)_{\lambda}$ for an arbitrary $a \in R$.
\end{lema}

{\bf Proof.} Denote $w = e_{w_4-w_1}(a)v_{\lambda} -
e_{w_4-w_1}(av)_{\lambda}$.  Clearly, $e_{w_1-w_4}(R)w = (0)$.
Since $f(w_3 - w_4) >0$, it follows that $e_{w_3-w_4}(b)w =
e_{w_3-w_1}(ab)v_{\lambda}  - e_{w_3-w_1}(b)(av)_{\lambda} = 0$ by
Lemma 6.4.

From $[q_{w_1+w_4}, e_{w_4-w_1}(R)] = (0)$ we conclude that
$q_{w_1+w_4}w = 0$.  Hence, $[q_{w_1+w_4}, e_{w_3-w_4}(R)]w =
(0).$ Also, $q_{-w_2-w_4}(R)w = (0)$ by Lemma 3.2 applied to the
root $-w_1-w_2$.  We proved that $L_{w_1+w_3}w = (0)$.

If $w  \neq 0$ then $U(L_+)w \cap V_{\lambda} \neq (0).$

It means that there exist positive roots $\alpha_1, \ldots,
\alpha_s$ such that $\alpha_1 + \cdots + \alpha_s + w_4 - w_1 \in
{\bf Z}(w_1+w_2+w_3+w_4)$ and, moreover, $\alpha_i+w_4 - w_1$ is a
negative root or 0 for any $i$. If $\alpha_i$ is an even root and
$\alpha_i+w_4 - w_1$  is one of the roots of Lemma 3.2 or 0 then
$e_{\alpha_i}(b) w = 0$ again by Lemma 3.2.  By Lemma 6.4
$\alpha_i$ is not supposed to be $w_3 - w_4$ as well. This rules
out all even roots except $w_2-w_4$.

Of odd roots, we have to examine $w_1+w_2$ and $w_1+w_3$, but the
latter one has been ruled out above.  Hence, $i(w_2-w_4) +
j(w_1+w_2) + (w_4 - w_3) = k(w_1 + w_2 + w_3 + w_4)$; $i,j,k \in
{\bf Z}$; $i,j \geq 0$.  This equation does not have a solution.
Hence $w = 0$.  This finishes the proof of the lemma.

\begin{lema}  $e_{w_3-w_1}^2 v_{\lambda} = 0$.
\end{lema}

{\bf Proof.} For an arbitrary element $a \in R$ we have
$e_{w_1-w_3}(a)e_{w_3-w_1}^2 v_{\lambda} = [e_{w_1-w_3}(a),
e_{w_3-w_1},e_{w_3-w_1}] v_{\lambda} + 2 e_{w_3-w_1}h_{w_1-w_3}(a)
v_{\lambda} = -2 e_{w_3 - w_1}(a)v_{\lambda} +$ \break $ 2
e_{w_3-w_1} h_{w_1-w_3}(a) v_{\lambda} = 0$ by Lemma 6.5.

\medskip

Now, as in Section 4 we see that $$U(L_+)e^2_{w_3 -
w_1}v_{\lambda} \cap V_{\lambda} = L^2_{w_1+w_2}
L^2_{w_1+w_4}e^2_{w_3-w_1} v_{\lambda}.$$  We have  $L_{w_1+w_4} =
[q_{w_1+w_3},e_{w_4 -w_3}(R)] + q_{-w_2-w_3}(R)$.

\medskip

{\bf Claim 1:} $[q_{w_1+w_3}, e_{w_4-w_3}(a)]e^2_{w_3-w_1}
v_{\lambda} = [q_{w_1+w_3},
e_{w_4-w_3}]e^2_{w_3-w_1}(av)_{\lambda}$ for an arbitrary $a \in
R$.

Indeed, $[L_{w_1+w_4},e_{w_3 - w_1},e_{w_3 - w_1}] = (0)$ implies
$$[q_{w_1+w_3}, e_{w_4 - w_3}(a)]e_{w_3-w_1}^2 v_{\lambda}  = 2e_{w_3 - w_1} [[q_{w_1+w_3}, e_{w_4-w_3}(a)],e_{w_3 -
w_1}] v_{\lambda}=$$
$$2 e_{w_3 - w_1}q_{w_1+w_3} e_{w_4-w_1}(a) v_{\lambda} = 2e_{w_3 - w_1} q_{w_1+w_3} e_{w_4-w_1}(av)_{\lambda}$$  by
Lemma 6.5.  This proves the claim.

\medskip

{\bf Claim 2:}  $q_{-w_2-w_3}(a)e^2_{w_3-w_1} v_{\lambda} =
q_{-w_2-w_3}e_{w_3-w_1}^2 (av)_{\lambda}$.

Indeed, $q_{-w_2-w_3}(a)e_{w_2-w_1}^2 v_{\lambda} = 2
e_{w_3-w_1}[q_{-w_2-w_3}(a),e_{w_3-w_1}]v_{\lambda}= $ \break
$2e_{w_3-w_1}q_{-w_1-w_2}(a) v_{\lambda} = 2
e_{w_3-w_1}q_{-w_1-w_2}(av)_{\lambda}$ by Lemma 3.2.  This proves
the claim.

These claims and the similar assertions for $L_{w_1+w_2}
v_{\lambda}$ show that $$L^2_{w_1+w_2} L_{w_1+w_4}^2 e^2_{w_3-w_1}
V_{\lambda} =$$  $$ (Fq_{w_1+w_2} + Fq_{-w_3-w_4})^2(Fq_{w_1+w_4}
+ F q_{-w_2-w_3})^2e^2_{w_3-w_1} V_{\lambda} =
$$ $$q_{w_1+w_2}q_{w_1+w_4} q_{-w_3-w_4}q_{-w_2-w_3} V_{\lambda}.$$

The computations of Lemma 6.3 show that, under the assumption
\break $<\lambda, h_{w_2-w_3}> = 0$, this expression is equal to
0. This finishes the proof of the lemma.

\begin{lema} $e_{w_2-w_3} v_{\lambda} = 0$.
\end{lema}

{\bf Proof.} If $e_{w_2-w_3} v_{\lambda} \neq 0$ then there exist
positive roots $\alpha_1, \ldots, \alpha_s$ such that  $\alpha_1 +
\cdots + \alpha_s = w_3-w_2$, for each $\alpha_1$ the sum
$\alpha_i + w_2-w_3$ is a negative root or 0 and $L_{\alpha_1}
\cdots L_{\alpha_s}e_{w_2-w_3}v_{\lambda} \neq (0).$

The only positive roots with the properties above are $w_3-w_2$
and $w_1 + w_4$.  But $$e_{w_3-w_2}(a)e_{w_2-w_3}v_{\lambda} =
h_{w_3-w_2}(a) v_{\lambda} = < \lambda, w_3-w_2>(av)_{\lambda} =
0.$$ Hence, all $\alpha_i$ have to be equal to $w_1+w_4$,
$L^s_{w_1+w_4}e_{w_2-w_3}v_{\lambda} \neq (0)$. But $[L_{w_1+w_4},
[L_{w_1+w_4}, e_{w_2-w_3}]] = (0)$, which leads to a contradiction
and finishes the proof of the lemma.

\begin{lema}  $e_{w_4 - w_3}(a)v_{\lambda} =
e_{w_4-w_3}(av)_{\lambda}.$
\end{lema}

We have $e_{w_4-w_3}(a) = [e_{w_4 - w_1}(a), e_{w_1-w_3}]$.  Hence
$$e_{w_4-w_3}(a) = - e_{w_1-w_3}e_{w_4-w_1}(a) v_{\lambda} =$$
$$-e_{w_1-w_3}e_{w_4-w_1}(av)_{\lambda} = [e_{w_4-w_1},
e_{w_1-w_3}](av)_{\lambda} = e_{w_4-w_3}(av)_{\lambda}$$ by Lemma
6.5.

\begin{lema}  $e_{w_4-w_3}^{<\lambda,w_3-w_4>+1} v_{\lambda} = 0$.
\end{lema}

{\bf Proof.}  As in the proof of Lemma 6.7, if the assertion is
not true, then there exist positive roots $\alpha_1, \ldots,
\alpha_s$, such that $\alpha_1 + \cdots + \alpha_s + (<\lambda,
w_3-w_4> + 1)(w_4-w_3)\in {\bf Z}(w_1+ \cdots + w_4)$, $\alpha_i +
w_4 - w_3$ is a negative root or 0.  In fact, 0 is also excluded,
because
$$e_{w_3-w_4}(a)e_{w_4-w_3}^{<\lambda,w_3-w_4>+1} v_{\lambda} =
e_{w_3-w_4}e_{w_4-w_3}^{<\lambda,w_3-w_4>+1} (av)_{\lambda} = 0$$
by Lemma 6.8.

\medskip

The only such positive roots are: $w_3 - w_2, \, w_2-w_4$,
$w_1+w_2$, $-2w_4$.

\medskip

Hence, there exist $i,j,k,l \in {\bf Z}_{\geq 0}$, $p \in {\bf
Z}$, such that $i(w_3-w_2) + j(w_2-w_4) + k(w_1+w_2) - 2lw_4 +
m(w_4-w_3) = p(w_1+w_2+w_3+w_4)$, where $m = <\lambda,
w_3-w_4>+1$.  It means, that $k = p$, $-i+j+k = p$, $i - m = p$,
$-j -2l = p$.  The first two equalities imply that $p = k \in {\bf
Z}_{\geq 0}$, $i = j$.  Now, adding the last two equalities we get
$-m - 2l = 2p$, where the left hand side is negative, whereas the
right hand side is positive.  This concludes the proof of the
lemma.

\bigskip

The element $e_{\pm(w_1-w_3)}$, $e_{\pm(w_2-w_3)}$,
$e_{\pm(w_4-w_3)}$ generate $sl(4)$ and act on $V_{\lambda}$
nilpotently.  Arguing as in the proof of Lemmas 4.11 and 4.12 we
get

\begin{lema} (1) There exists $m \geq 1$ such that $e_{w_i-w_j}^m
V_{\lambda} = (0)$ for any $1 \leq i \neq j \leq 4$.

(2) There exist $m \geq 1$ such that $L_{\alpha}^m V_{\lambda} =
(0)$ for an arbitrary root $\alpha$.

\end{lema}

This implies that for an integral dominant weight $\lambda$ such
that $ 1 = \break <\lambda, h_{w_1-w_3}> $, $0 = <\lambda,
h_{w_3-w_2}>$, the module $V(\lambda, -1,\alpha)$ is a finite sum
of weights spaces with respect to the action of $H$.

\section{ The case $<\lambda, h_{w_1-w_3}> = 0$}

We have $$[[e_{w_4-w_1}(a),q_{w_3+w_1}],q_{w_2+w_1}] =
-[[e_{w_3-w_1}(a), q_{w_1+w_4}],q_{w_2+w_1}] = Vir(a).$$

If $<\lambda, h_{w_1-w_3}> = 0$ and nevertheless $V(\lambda,
\beta,\alpha)$ is of finite type, then $e_{w_3-w_1} V_{\lambda} =
(0)$.  This implies that $$e_{w_3-w_1}(a) V_{\lambda} = - {1 \over
2} [e_{w_1-w_3}(a), e_{w_3-w_1},e_{w_3-w_1}] V_{\lambda} = (0).$$
Hence \hskip 1 true cm  $$q_{w_2+w_1}q_{w_1+w_4}e_{w_3-w_1}(a)
V_{\lambda} = [q_{w_2+w_1},[q_{w_1+w_4}, e_{w_3-w_1}(a)]]
V_{\lambda} = Vir(a) V_{\lambda} = (0).$$ Since $H \subseteq [H
\otimes R, Vir(R)]$ it follows that $HV_{\lambda} = (0)$, $\lambda
= 0$.  Then $V$ is a 1-dimensional module with zero
multiplication, which is not viewed as irreducible.

This concludes the proof of Theorem 3.1.

\section{Jordan bimodules}

Let $V$ be a Jordan bimodule over a unital Jordan (super)algebra
$J$.  Then $V$ can be represented as a direct sum $V = V_0 \oplus
V_{1/2} \oplus V_1$, where $JV_0 = (0)$, $V_{1/2}$ is a one-sided
Jordan bimodule (see [MZ2] ), $V_1$ is a unital Jordan bimodule.

\medskip

In [MZ2] it was shown that the universal associative enveloping
algebra of $JCK(R,d)$ is $M_4(W(R,d)$.  It means that one-sided
Jordan $JCK(R,d)$-bimodules are left modules over $M_4(W(R,d))$
or, equivalently, $4$-tuples $U^4$, where $U$ is a left module
over $W(R,d)$.

\medskip

Now suppose that $R = F[t^{-1},t]$, $d = {d \over dt}$, and the
one-sided Jordan bimodule over $JCK(6)$ is conformal.Then the left
module $M$ over $W(R,d)$ is a unital conformal module. Such
modules correspond to left unital ${\bf C}[d]$-modules.
Irreducible ${\bf C}[d]$-modules are one-dimensional and
parametrized by scalars $\alpha \in F$. Indecomposable conformal
modules of finite type correspond to Jordan blocks.

\medskip

Let's be more precise.  let $N$ be a left ${\bf C}[d]$-module.
Then $N[t^{-1},t] = \{\sum n_it^{i}, n_i \in N, i \in {\bf Z} \}$
is a conformal left $W(F[t^{-1},t], {d \over dt})$-module.  It is
of finite type if and only if $dim_{\bf C} N < \infty$.  The space
of 4-tuples $N[t^{-1},t]^4$ is a left associative conformal module
over $M_4(W)$, hence a one-sided Jordan conformal module over
$JCK(6)$.

\begin{prop}  (1) Every one-sided Jordan conformal bimodule over
$JCK(6)$ is of the type $N[t^{-1},t]^4$;

(2) The module $N[t^{-1},t]^4$ is irreducible if and only if $N$
is one-dimensional, $N = {\bf C}v$, $vd = \alpha v$, $\alpha \in
F$;

(3) $N[t^{-1},t]^4$ is an indecomposable module of finite type if
and only if $N$ is a Jordan block.
\end{prop}

\bigskip

Now let $V$ be a unital irreducible conformal Jordan bimodule of
finite type over $J = JCK(6)$.  In [MZ4], [Z] it was shown that
the Tits-Kantor-Koecher construction $K(V)$ is a Lie conformal
module of finite type over the  TKK-algebra $K(J) = CK(6)$.  In
[MZ4] we proved that: (1) the reduced module $\bar K(V)$ is
irreducible over $K(J)$ and uniquely determines the $J$-bimodule
$V$; (2) the action of the Cartan subalgebra $H$ on $\bar K(V)$ is
diagonalizable, all weights of $\bar K(V)$ belong to the set
$\{\pm w_i \pm w_j, 1 \leq i,j \leq 4 \}$.  Let $\lambda$ be the
highest weight of $\bar K(V)$.  Since the Weyl group of $sl_4$ is
the permutation group $P_4$ and $f(\lambda) \geq
f(\sigma(\lambda))$ for all $\sigma \in P_4$, the only
possibilities for $\lambda$ are : $2w_1$, $w_1-w_4$, $w_1+w_3$,
$-2w_4$.  The last two cases are ruled out by Theorem 3.1. The
modules $V(2w_1, \beta,\alpha)$, $\beta, \alpha \in F$; and
$V(w_1-w_4, -1, \alpha)$, $\alpha \in F$ indeed have Jordan
structures.

\begin{prop}  Unital irreducible conformal Jordan
$JCK(6)$-bimodules of finite type form two parametric families,
which correspond to $V(2w_1, \beta,\alpha)$, $\beta,\alpha \in F$
and $V(w_1-w_4, -1,\alpha)$, $\alpha \in F$.
\end{prop}

{\bf Proof.} We need to show that $V(2w_1, \beta, \alpha)$,
$V(w_1-w_4, -1, \alpha)$, $\alpha,\beta \in F$ are reduced
Tits-Kantor-Koecher modules of the form $\bar K(V)$ for some
unital Jordan bimodules $V$ over $J = JCK(6)$.

\medskip

As in section 5 consider the associative commutative algebra
$\tilde R  = R + Rv$, $R = {\bf C}[t^{-1},t]$, $v^2 = 0$ with the
derivation $d$, $d(t) = -1$, $d(v) = \alpha v$.  Consider the
algebra $\tilde W(\tilde R, d) = \sum_{\beta \in F}
W_{\beta}(\tilde R, d)$, $W_{\beta}(R,v,d) = \sum_{i = 0}^{\infty}
Rvd^{\beta - i}$.  The Lie superalgebra $L = CK(6)$ is embedded
into $M_8(W)$, hence into $M_8(\tilde W)$.

\medskip

Consider the subspace $(W_{\beta}(R,v,d))_{1,5}$ of matrices
having $W_{\beta}(R,v,d)$ at the intersection of the 1st row  and
5th column and 0 elsewhere. It is easy to see that $[L_+,
(W_{\beta}(R,v,d))_{1,5}] = (0)$ and for an arbitrary element  $u
\in (W_{\beta}(R,v,d))_{1,5}$, arbitrary $1 \leq i \neq j \leq 4$,
we have $[h_{w_i-w_j}, u] = (2w_1 / w_i-w_j)u$. Let $U_{\beta} =
U(L)(W_{\beta}(R,v,d))_{1,5}$ be the $L$-submodule generated by
$(W_{\beta}(R,v,d))_{1,5}$ in $M_8(\tilde W)$.  Clearly, $2w_1$ is
the highest weight of this submodule and $(U_{\beta})_{2w_1} =
(W_{\beta}(R,v,d))_{1,5}$.

\medskip

The $L$-submodule $U_{\beta - 1}$ of $U_{\beta}$ is generated by
$(W_{\beta - 1}(R,v,d))_{1,5}$, $U_{\beta} / U_{\beta - 1} = U(L)
(W_{\beta}(R,v,d)/ W_{\beta - 1}(R,v,d))$, hence $U_{\beta} /
U_{\beta - 1} = U(L)(U_{\beta} / U_{\beta - 1})_{2w_1}$, and $
(U_{\beta} / U_{\beta - 1})_{2w_1} \simeq V(\beta,\alpha)$.

Hence $V(2w_1,\beta,\alpha)$ is a homomorphic image of the module
$U_{\beta} / U_{\beta - 1}$.  Then all weights of
$V(2w_1,\beta,\alpha)$ belong to the set $\{\pm w_i \pm w_j. 1
\leq i,j \leq 4 \}$, which implies that $V(2w_1,\beta,\alpha)
\simeq K(\bar V)$ for some irreducible unital Jordan $J$-bimodule
$V$.

\medskip

Now let's turn to bimodules $V(2w_1,\beta,\alpha)$. Consider the
Cheng-Kac superalgebra $CK(\tilde R, d)$ and the subspace
$e_{w_1-w_4}(Rv)$.  This subspace generates the $L$-submodule
which is isomorphic to $V(w_1-w_4, -1, \alpha)$.  This concludes
the proof of the proposition.

 \vskip 1cm

\centerline{\Large \bf References}

\vskip 0.5cm

[BKL1] C. Boyallian, V. G. Kac and J. I. Liberati, Irreducible modules
over finite simple Lie conformal superalgebras of type K , {\em J. Math. Phys.} {\bf 51},063507, (2010).

\medskip
[BKL2] C. Boyallian, V. G. Kac and J. I. Liberati, Classification of finite irreducible modules
over the Lie conformal superalgebra $CK_6$, Preprint.

\medskip

[BKLR] C. Boyallian, V. G. Kac,  J. I. Liberati and A. Rudakov, Representations of simple finite Lie conformal superalgebras
of type W and S, {\em J. Math. Phys.} {\bf 47}, 043513, (2006).

\medskip

[CK1] S. Cheng and V.G. Kac, Conformal modules, {\em Asian J. Math} {\bf 1}, 181-193, (1997).
Erratum: 2 ,153-156, (1998).

\medskip

[CK2] S. Cheng and V.G. Kac, A new N = 6 superconformal algebra, {\em Comm. Math. Phys.} {\bf 186}
no. 1, 219-231, (1997).

\medskip

[FK] D. Fattori and V.G. Kac, Classivication of finite simple Lie
conformal superalgebras, {\em J. of Algebra} {\bf
258}, no.1, 23-59,(2002).

\medskip

[GLS] P. Grozman, D. Leites and I. Shchepochkina, Lie superalgebras of
string theories, {\em Acta Math. Vietnam} {\bf 26} no. 1, 27-63, (2001).

\medskip

[K1] V.G. Kac, Classification of simple Z-graded Lie superalgebras and
simple Jordan superalgebras, Comm. in Algebra 5 (13), 1375-1400, (1977).

\medskip

[K2] V.G. Kac, Vertex algebras for beginners, University Lecture Series 10, American Mathematical Society, Providence, RI
1996, Second Edition  1998.

\medskip

[MZ1] C. Mart\'{\i}nez and E. Zelmanov, Simple finite-dimensional Jordan
Superalgebras in prime characteristic , {\em Journal of Algebra} {\bf 236} no.2, 575-629,
(2001).

\medskip

[MZ2] C. Mart\'inez and E. Zelmanov, Specializations of Simple Jordan
Superalgebras, {\em Canad. Math. Bull.}  {\bf 45} (4), 653-671, (2002).

\medskip

[MZ3] C. Mart\'inez and E. Zelmanov, Lie superalgebras graded by P(n)
and Q(n), {\em Proc. Natl. Acad. Sci. USA}  {\bf 100} no.14, 8130-8137, (2003).

\medskip

[MZ4]C. Mart\'inez and E. Zelmanov, Representation Theory of Jordan Superalgebras I, {\em Trans. Amer. Math. Soc} {\bf 362} no. 2, 815-846, (2010).

\medskip

[Z] E. Zelmanov, On the structure of conformal algebras, {\em Contemporary Math.} {\bf 264}, 139-153 , (2000).

\end{document}